\documentclass[english]{amsart}

\usepackage{amsmath,amssymb,enumerate,bbm}

\usepackage[T1]{fontenc}
\usepackage[utf8]{inputenc}
\usepackage[all]{xy}

\usepackage{babel}
\usepackage{amstext}
\usepackage{amsmath}
\usepackage{amsfonts}
\usepackage{latexsym}
\usepackage{ifthen}

\usepackage{xypic}
\xyoption{all}
\newcommand{\cal}{\mathcal}

\newcommand{\PN}{{\mathbb P}}

\newcommand{\rk}{{\rm rk}}

\newtheorem{lemma1}{}[section]

\newenvironment{lemma}{\begin{lemma1}{\bf Lemma.}}{\end{lemma1}}
\newenvironment{example}{\begin{lemma1}{\bf Example.}\rm}{\end{lemma1}}

\newenvironment{theorem}{\begin{lemma1}{\bf Theorem.}}{\end{lemma1}}
\newenvironment{proposition}{\begin{lemma1}{\bf Proposition.}}{\end{lemma1}}
\newenvironment{corollary}{\begin{lemma1}{\bf Corollary.}}{\end{lemma1}}
\newenvironment{remark}{\begin{lemma1}{\bf Remark.}\rm}{\end{lemma1}}
\newenvironment{remarks}{\begin{lemma1}{\bf Remarks.}\rm}{\end{lemma1}}
\newenvironment{definition}{\begin{lemma1}{\bf Definition.}}{\end{lemma1}}
\newenvironment{notation}{\begin{lemma1}{\bf Notation.}}{\end{lemma1}}

\newenvironment{conjecture}{\begin {lemma1}{\bf Conjecture.}}{\end{lemma1}}

\newenvironment{remark*}{{\bf Remark.}}{}
\newenvironment{remarks*}{{\bf Remarks.}}{}
\newenvironment{example*}{{\bf Example.}}{}

\newcommand{\Q}{\ensuremath{\mathbb{Q}}}
\newcommand{\Z}{\ensuremath{\mathbb{Z}}}
\newcommand{\C}{\ensuremath{\mathbb{C}}}
\newcommand{\N}{\ensuremath{\mathbb{N}}}
\newcommand{\PP}{\ensuremath{\mathbb{P}}}
\newcommand{\D}{\ensuremath{\mathbb{D}}}

\newcommand{\merom}[3]{\ensuremath{#1:#2 \dashrightarrow #3}}

\newcommand{\holom}[3]{\ensuremath{#1:#2  \rightarrow #3}}
\newcommand{\fibre}[2]{\ensuremath{#1^{-1} (#2)}}

\makeatletter
\ifnum\@ptsize=0  \fi \ifnum\@ptsize=2  \fi 

\newcommand\sO{{\mathcal O}}

\DeclareMathOperator*{\red}{red}

\newcommand{\chow}[1]{\ensuremath{\mathcal{C}(#1)}}

\title{Non-algebraic  compact K\"ahler threefolds \\ admitting endomorphisms} 
\date{July 21, 2009}

\author{Andreas H\"oring}
\author{Thomas Peternell}

\address{Andreas H\"oring, Universit\'e Paris 6, Institut de Math\'ematiques de Jussieu, Equipe de Topologie et G\'eom\'etrie Alg\'ebrique, 175, rue du Chevaleret, 75013 Paris, France}
\email{hoering@math.jussieu.fr}

\address{Thomas Peternell, Mathematisches Institut, Universit\"at Bayreuth, 95440 Bayreuth, 
Germany 
}
\email{thomas.peternell@uni-bayreuth.de}

\begin{document}

\begin{abstract}
We classify non-algebraic  compact K\"ahler threefolds admitting an endomorphism 
$f: X \rightarrow X$ of degree at least two.
\end{abstract}

\subjclass[2000]{32J17, 32H02, 32J27, 14E20, 14E30}
\keywords{endomorphism, compact K\"ahler manifold, non-algebraic manifold, torus fibrations, MMP}

\maketitle

\vspace{-1cm}

\tableofcontents

\vspace{-1cm}

\section{Introduction}

\medskip

An endomorphism of a compact complex manifold $X$ is a surjective map \holom{f}{X}{X},
usually assumed to be of degree $d$ at least two, i.e.\ automorphisms are excluded.
Endomorphisms of projective manifolds were intensively studied in the last years \cite{Bea01, Fuj02, Nak02, HM03, Ame03, FN05, FN07, NZ07, Nak08, AKP08, NZ09}.
For example if $X$ is a smooth projective threefold and $f$ is \'etale, $X$ is completely classified up to \'etale cover. Also the higher-dimensional case
is intensively treated but far from being completely understood.

In this paper we classify all non-algebraic three-dimensional compact K\"ahler manifolds $X$ admitting an endomorphism $f$ regardless
whether $f$ is ramified or not.
Before we state our classification results let us give an example 
how the non-algebraicity assumption gives additional restrictions
on the existence of endomorphisms. 

\begin{theorem} \label{theoremkappanminusone}
Let $X$ be a non-algebraic compact K\"ahler manifold of dimension $n$ which admits a meromorphic endomorphism \merom{f}{X}{X} of degree $d>1$.
Then $\kappa(X) \ne n-1$.
\end{theorem}

In fact if $\kappa (X) = n-1$, the Iitaka fibration coincides with the algebraic reduction, which is an elliptic fiber space over a projective manifold
of dimension $n-1$. The existence of the endomorphism gives rise to a meromorphic multi-section of the fibration, and therefore any two points of $X$ can be joined by a chain of compact curves. Hence $X$ is projective due to a result of Campana.

If the manifold is not uniruled our results can be summarised as follows.

\begin{theorem} \label{maintheoremA}
Let $X$ be a compact non-algebraic K\"ahler threefold which is not uniruled. 
Suppose that $X$ admits an endomorphism $f: X \to X$ of degree $d > 1.$ 
Then (up to \'etale cover) one of the following holds:
\begin{enumerate}
\item $\kappa(X)=0:$ then either
\begin{enumerate}
\item $X$ is a torus or 
\item $X$ is a product $Y \times E$ where $Y$ is bimeromorphic to a torus or K3 surface and $E$ an elliptic curve.   
\end{enumerate}
\item $\kappa(X)=1:$ then either
\begin{enumerate}
\item
$X$ is a product $C \times A$ where $C$  is a curve of general type and $A$ is a two-dimensional torus of algebraic dimension at most one or 
\item $X$ is a product $E \times S$ where $E$  is an elliptic curve and $S$ 
an elliptic surface of algebraic dimension one.
\end{enumerate}
\end{enumerate}
\end{theorem}

If $f$  is ramified, it is a basic fact that $K_X$ is not pseudoeffective. If $X$ is algebraic, then $X$ is uniruled (cf. \cite{BDPP04} which is based on a theorem of Miyaoka-Mori using characteristic $p$ methods).
In the K\"ahler case the uniruledness is only known in dimension three due to a remarkable result of Brunella \cite{Bru06}.

\begin{theorem} \label{maintheoremB}
Let $X$ be a compact non-algebraic K\"ahler threefold of algebraic dimension $a(X)$ which is uniruled. 
Suppose that $X$ admits an endomorphism $f: X \to X$ of degree $d > 1.$
If $f$ is ramified, then (up to \'etale cover) one of the following holds.
\begin{enumerate}
\item $a(X)=0:$ then $X$ is a projectivised bundle $\PP(E)$ over a torus $A$ of algebraic dimension zero,
$f$ induces an endomorphism on $A$ of degree at least two, and $E$ is a direct sum of line bundles.
\item $a(X)=1:$ then either
\begin{enumerate}
\item $X$ is a product $S \times \PP^1$, where $S$ is  a compact K\"ahler surface of algebraic dimension zero
and $f$ induces an automorphism on $S$ or
\item $X$ is a projectivised bundle $\PP(E)$ over a torus $A$ of algebraic dimension at most one,
$f$ induces an endomorphism on $A$ of degree at least two, and $E$ is a direct sum of line bundles.
\end{enumerate}
\item $a(X)=2:$ then either
\begin{enumerate}
\item $X$ is a product $Y \times \PP^1$ where $Y$ is a surface of algebraic dimension one and 
$f$ induces an automorphism on $Y$ or
\item $X$ is a projectivised bundle $\PP(E)$ over a torus $A$ of algebraic dimension one and $f$ induces
an endomorphism $g$ of degree at least two on $A$.
\end{enumerate}
\end{enumerate}
If $f$ is \'etale, then $X$ is (up to \'etale cover) a 
projectivised bundle $\PP(E)$ over a non-algebraic torus $A$ and $c_1^2(E)=4c_2(E)$.
\end{theorem}

In the projective case Mori theory is heavily used to pass to minimal models. 
In the K\"ahler case only rudiments of Mori theory are known (see Section \ref{subsectionmori}), but they suffice
in the special cases we are interested in. Moreover the algebraic reduction provides a powerful tool.

{\em Guide to the reader.}
In Section \ref{sectionbasic} we recall the basic definitions and 
gather various results which will be used at some point in the paper.  In Section \ref{sectionthreefolds}
we prove new statements in the Mori theory of compact K\"ahler threefolds which should be interesting
in their own right. Together with
the results from \cite{Pet98,Pet01} they allow us to establish a MMP for compact K\"ahler threefolds
admitting \'etale endomorphisms of degree at least two. A crucial point is that the contractions never contract a
divisor to a point, so that we always stay in the smooth category. 
Section \ref{sectiontorusfibrations} provides results that would be trivial in the projective case:
based on a discussion of the fixed point set of torus endomorphisms, we establish the existence
of multisections for torus fibrations commuting with an endomorphism. Using a theorem of Nakayama and Zhang, the proof of 
Theorem \ref{theoremkappanminusone}
is then an easy exercise (cf. page \pageref{prooftheoremkappanminusone}). The Theorems \ref{maintheoremA} and \ref{maintheoremB} are proven in
the Sections \ref{sectionnonuniruled} and \ref{sectionuniruled} respectively. 
These two sections are the core of this paper, we advise the reader to start here
and skip the preceding technical sections for the first reading.

{\bf Acknowledgements.} We would like to thank T.-C. Dinh, C. Favre, P. Popescu-Pampu and N. Sibony
for very helpful discussions on this subject. We also thank the Research Group ``Classification of Algebraic Surfaces and Compact Complex Manifolds''
of the Deutsche Forschungsgemeinschaft DFG for financial support.

\section{Notation and basic results} \label{sectionbasic}

For standard definitions in complex geometry 
we refer to  \cite{Ha77} or \cite{Kau83}.  Moreover we refer to \cite{BHPV04} for basic results on surfaces
and to \cite{Fuj83} to the classification theory of higher-dimensional non-algebraic varieties. 
Manifolds and varieties are always supposed to be irreducible.
We will always assume implicitly that a compact K\"ahler manifold/surface/threefold is smooth.
If a certain statement holds for a singular variety, we will mention what types of singularities are allowed.

We  say that a certain property holds for a general (resp. very general) point $x \in X$ if
there exists a finite (resp. countable) union of proper subvarieties of $X$ such that the property 
holds for every point in the complement.

\subsection{Endomorphisms}

\begin{notation} Let $X$ be a compact complex variety that is normal or Gorenstein. An endomorphism
is a holomorphic surjective map $f: X \rightarrow X$. 
It is easy to see that $f$ is a finite map, so the 
ramification formula 
\[
K_X = f^* K_X + R
\]
holds and we will call the support of $R$ the ramification locus.
The support of the cycle theoretic image $B:=f_* R$ will be called the branch locus.
We will say that $f$ is \'etale (in codimension one) if $R$ is empty.
\end{notation}         

\begin{remarks*} 1. If $X$ is smooth, an endomorphism that is \'etale in codimension one is \'etale in every point.

2. Analogously one defines a meromorphic endomorphism of a compact complex variety as 
a meromorphic dominant map $f: X \dashrightarrow X$. 
\end{remarks*}

The following well-known results will be used many times in this paper:

\begin{lemma} 
Let $X$ be a compact K\"ahler manifold of dimension $n$, and let \holom{f}{X}{X} be an endomorphism of degree $d>1$. 
Then $f$ is finite. The linear maps
$$
\holom{f^*}{H^*(X, \Q)}{H^*(X, \Q)} \ \mbox{and} \ \holom{f_*}{H^*(X, \Q)}{H^*(X, \Q)}
$$
are isomorphisms. More precisely we have $f_*  f^* = d \ {\rm Id}$.

If $f$ is \'etale, we have
\begin{itemize}
\item $\chi(X, \sO_X)=0$,
\item $e(X) := \chi_{\rm top}(X) = c_n(X) = 0$, and
\item $K_X^n=0$.
\end{itemize}
\end{lemma}

\begin{lemma} \label{lemmaintersection}
Let $X$ be a compact K\"ahler variety of dimension $n$, and let \holom{f}{X}{X} be an endomorphism of degree $d>1$. Let $D$ be a Cartier divisor on $X$ such that
$f^* D \equiv_{num} m D$ for some $m \in \N$. Then we have $D^n=0$ or $d=m^n$.

In particular if $D$ is an effective divisor that is not contained in the branch locus of $f$
and such that $\fibre{f}{D}=D$, we have $D^n=0$.
\end{lemma}

\begin{proof}
Since $f_* f^* = d \ Id$, we see that
$$
m^n D^n = (f^* D)^n = f^*(D^n) = d D^n,
$$
so the first statement is immediate. For the second statement observe that the hypothesis implies that $f^* D \simeq D$.
\end{proof}

The next statement shows that from the point of view of endomorphisms, it is quite natural to treat separately uniruled and non-uniruled manifolds.

\begin{proposition} \label{propositionuniruled} Let $X$ be a smooth compact K\"ahler threefold and $f: X \to X$ be a ramified endomorphism of degree $d>1$. 
Then $X$ is uniruled.
\end{proposition} 

\begin{proof} 
By \cite[Thm.4.1]{AKP08} the canonical divisor $K_X$ is not pseudo-effective, 
i.e.\ the class of $K_X$ is not contained in the closure of the K\"ahler cone. Therefore by \cite[Cor.1.2]{Bru06},
$X$ is uniruled.
\end{proof}   

\subsection{Fibrations}

A fibration is a proper surjective morphism \holom{\varphi}{X}{Y} with connected fibres
from a complex manifold onto a normal complex variety $Y$. 
A fibre is always a fibre in the scheme-theoretic sense and will be denote by $\fibre{\varphi}{y}$ or $X_y$.
A set-theoretic fibre is the reduction of the fibre.
The $\varphi$-smooth locus is the largest Zariski open subset $Y^* \subset Y$
such that for every $y \in Y^*$, the fibre $\fibre{\varphi}{y}$ is a smooth variety of dimension $\dim X - \dim Y$.
The $\varphi$-singular locus is its complement.

Let us recall the rigidity lemma that will be used many
times in this paper.

\begin{lemma} \label{lemmarigidity}
Let $\holom{f}{X}{Y}$ and $\holom{g}{X}{Z}$ be fibrations. Suppose that for every $z \in Z$ the fibre $\fibre{g}{z}$
is mapped by $f$ onto a point. Then there exists a holomorphic map \holom{h}{Z}{Y} such that $f=h \circ g$.

If moreover $g$ is flat the same conclusion holds if at least one $g$-fibre is contracted by $f$.
\end{lemma}

\begin{definition}
\label{definitionalmostholomorphic}
A meromorphic map \merom{\varphi}{X}{Y} from a compact K\"ahler manifold to a normal K\"ahler variety
is almost holomorphic if  
there exist non-empty Zariski open subsets $X^* \subset X$ and $Y^* \subset Y$ 
such that \holom{\varphi|_{X^*}}{X^*}{Y^*} is a fibration.
In particular for $y \in Y$ a general point, the fibre \fibre{\varphi}{y} exists in the usual sense
and is compact.
\end{definition}

The importance of almost holomorphic maps is due to the fact that every compact K\"ahler manifold admits such a fibration that separates the rationally connected part and the non-uniruled part: the {\it rationally connected quotient}.

\begin{theorem} \cite[Thm.5.4]{Ko96},\cite[Thm. 1.1]{Ca04}, \cite{GHS03}\footnote{The statement in \cite{Ko96} is in the algebraic setting, but the same proof goes through in the compact K\"ahler category: the main technical tool \cite[Thm. 1.1]{Ca04} holds in this larger generality.}
\label{definitionrationallyconnectedquotient}
Let $X$ be a compact K\"ahler manifold. Then there exists an almost holomorphic fibration 
\merom{\varphi}{X}{Y} onto a normal compact K\"ahler variety $Y$ 
such that the general fibre is rationally connected and the variety $Y$ is 
not uniruled. 
This map is unique up to meromorphic equivalence of fibrations \cite{Ca04}
and will be called the rationally connected quotient.
\end{theorem}

\subsection{Endomorphisms that preserve fibrations}

\begin{definition}
Let $X$ be a compact K\"ahler manifold, and let \holom{f}{X}{X} be an endomorphism of degree $d>1$.
Suppose that $X$ admits a fibration \holom{\varphi}{X}{Y} onto a normal K\"ahler variety $Y$. 
If there exists an endomorphism \holom{g}{Y}{Y} such that $g \circ \varphi=\varphi \circ f$, 
we say that $f$ preserves the fibration and $g$ is the endomorphism induced by $f$ on $Y$.
\end{definition}

All the natural fibrations attached to a variety are preserved by an endomorphism.

\begin{proposition}  \label{propositioncommutation}
Let $X$ be a compact K\"ahler manifold, and let \holom{f}{X}{X} be an endomorphism of degree $d>1$.

1. Let \holom{\alpha}{X}{Alb(X)} be the Albanese map.
Then there exists an endomorphism \holom{g}{Alb(X)}{Alb(X)} such that $g \circ \varphi=\varphi \circ f$.

2. Let \merom{\varphi}{X}{Y} be the Iitaka fibration.
Then there exists a automorphism \holom{g}{Y}{Y} such that $g \circ \varphi=\varphi \circ f$.

3. Let \merom{\varphi}{X}{Y} be the rationally connected quotient of $X$.
Then there exists a meromorphic map \merom{g}{Y}{Y} such that $g \circ \varphi=\varphi \circ f$.
\end{proposition}

\begin{proof}
The first statement follows immediately from the universal property of the Albanese map, the second statement
can be shown as in \cite[Prop.2.5]{Fuj02} where the projectiveness assumption is actually  not used.
For the last statement, note that
up to replacing $Y$ by a Zariski open dense subset, we may suppose that the almost holomorphic map $\varphi$ is holomorphic.
By the rigidity lemma it is sufficient to show that a general $\varphi$-fibre $F$ is contracted by $\varphi \circ f$. 
Since $F$ is rationally connected and $Y$ is not uniruled, this is trivially true.
\end{proof}

Before we can prove an analogue of Proposition \ref{propositioncommutation}
for the algebraic reduction, we need one more definition.

\begin{definition} \label{definitionpolar}
Let $X$ be a compact K\"ahler manifold. An integral effective divisor $D \subset X$ is polar if there exists
a meromorphic function $\psi$ on $X$ such that $D \subset div(\psi)$. 
\end{definition}

\begin{remark}
If $X$ is projective, every divisor is polar. If the algebraic dimension is zero, every divisor is non-polar.
\end{remark}

\begin{proposition} \label{propositionalgebraicreduction}
Let $X$ be a compact complex manifold 
that admits a meromorphic endomorphism \merom{f}{X}{X}. 
Denote by  \merom{\varphi}{X}{Y} the algebraic reduction of $X$.
Then there exists a meromorphic endomorphism
$\holom{g}{Y}{Y}$ such that $\varphi \circ f = g \circ \varphi$.
\end{proposition}

\begin{remark*}
If $X$ has algebraic dimension one the variety $Y$ is a smooth compact curve, so the meromorphic
endomorphism extends to a holomorphic map.
\end{remark*}

\begin{proof}
Up to replacing $X$ by some bimeromorphic model we may suppose that the algebraic reduction is holomorphic.
The endomorphism $f$ acts by pull-back on the meromorphic function field
$$
\C(X) = \C(Y),
$$
and we define $g$ to be the meromorphic map corresponding to $f^*: \C(Y) \rightarrow \C(Y)$.
Since every polar divisor on $X$ is contained in a pull-back from $Y$ and the pull-back of a 
polar divisor is polar, one sees easily that $\varphi \circ f = g \circ \varphi$.
\end{proof}

\begin{definition}
Let $Y$  be a normal compact K\"ahler variety and let \holom{g}{Y}{Y} be an endomorphism
of degree at least two.
We say that $g$ is totally ramified in a point $y \in Y$, if the set-theoretical fibre $(\fibre{g}{y})_{\red}$  is a singleton. 
\end{definition}

\begin{proposition} \label{propositionspecialfibres} 
Let $X$ be a compact K\"ahler manifold, and let \holom{f}{X}{X} be an endomorphism.
Let \holom{\varphi}{X}{Y} be a surjective morphism onto a normal K\"ahler variety such that 
there exists an endomorphism \holom{g}{Y}{Y} such that $g \circ \varphi=\varphi \circ f$.

1. For $m \in \N$, set 
$$
T_m  := \{ y \in Y \ | \ \dim \fibre{\varphi}{y} > m-1 \}.
$$ 
Then we have (set-theoretically) $\fibre{g}{T_m}=T_m$. If $T_m$ is finite
and $g$ of degree at least two, then $g$ is totally ramified in every point of $T_m$.  

2. Set 
$$
R  := \{ y \in Y \ | \ \fibre{\varphi}{y}  \ \mbox{is reducible} \}.
$$ 
Then we have (set-theoretically) $\fibre{g}{R}=R$. If $R$ is finite
and $g$ of degree at least two, then $g$ is totally ramified in every point of $R$.  

3. Let
$$
\Delta  := \{ y \in Y \ | \ \fibre{\varphi}{y}  \ \mbox{is singular} \}
$$ 
be the $\varphi$-singular locus and denote by $R$ the branch locus of $g$. 
Then $\fibre{g}{\Delta}$ is contained set-theoretically in the union of $\Delta$ and the branch locus $R$.
If $g$ is \'etale, then (set-theoretically) $\fibre{g}{\Delta}=\Delta$.

In particular if $g$ is an \'etale map of degree at least two, then $T_m$, $R$ and $\Delta$ are either of positive dimension or empty.
\end{proposition}

\begin{proof}
1. Since $\fibre{g}{T_m}$ has at least as many irreducible components as $T_m$, it is sufficient to show that $\fibre{g}{T_m} \subset T_m$.
Yet if $t \in T_m$ and $y \in \fibre{g}{t}$, then 
$\fibre{\varphi}{y}$ surjects via $f$ on  $\fibre{\varphi}{t}$, so $\dim \fibre{\varphi}{y} \geq \dim \fibre{\varphi}{t} > m-1$.
If $T_m$ is finite, then $\# \fibre{g}{T_m} \geq \# T_m$ and equality holds if and only if $g$ is totally ramified in every point of $T_m$.

2. Analogous to 1.

3.  If $t \in \Delta$ and $y \in \fibre{g}{t}$ such that $y \not\in R$, then for every $x \in \fibre{\varphi}{y}$ we have
$$
\rk T_{\varphi,x} = \rk T_{g \circ \varphi, x} = \rk T_{\varphi \circ f, x} \leq \rk T_{\varphi, f(x)}.
$$
Since $\fibre{\varphi}{t}$ is singular, there exists a point $f(x) \in \fibre{\varphi}{t}$ such that the tangent map does not have maximal rank.
It follows that $y \in \Delta$.
\end{proof}

\begin{lemma} \label{lemmafibrationcurve}
Let $X$ be a compact K\"ahler variety, and let \holom{f}{X}{X} be an endomorphism of degree $d>1$.
Suppose that there exists a fibration $\holom{\varphi}{X}{C}$ onto a smooth curve $C$ and an  
automorphism \holom{g}{C}{C} such that $g \circ \varphi=\varphi \circ f$.
If $B_i$ is an irreducible component of the branch locus of $f$, then $B_i$ surjects onto $C$.
\end{lemma}

\begin{proof}
By Proposition \ref{propositionspecialfibres} the $\varphi$-singular locus $\Delta$ satisfies $\fibre{g}{\Delta}=\Delta$.
Since $\Delta$ is finite, we can suppose, up to replacing $f$ by $f^k$ that $g$ is the identity on $\Delta$. 
We argue by contradiction and suppose that  $\varphi(B_i)=c$, then it is an irreducible component of the fibre $X_c \simeq \varphi^* c$.
Suppose first that $X_c$ is not a smooth fibre, then $c \in \Delta$ which implies $g^*c=c$. Therefore we have
$$
X_c \simeq \varphi^* g^* c  \simeq f^* \varphi^* c  \simeq f^* X_c,
$$ 
in particular $m_i B_i \simeq f^* B_i \simeq B_i$. 
Since the ramification index $m_i$  is strictly larger than one, this implies that $B_i$ is homologous to zero, a contradiction.
If $X_c$ is a smooth fibre, then $X_c=B_i$ and $g^* c=c'$ with $c' \not\in \Delta$. The same computation shows that $X_{c'} \simeq f^* B_i=m_i B_i$.
Since $X_{c'}$ is smooth, so reduced, we get a contradiction.
\end{proof}

The following statement generalises \cite[Thm.1]{Ame03}, its proof is a mere adaptation of Amerik's proof
to our context.

\begin{proposition} \label{propositionamerik}
Let $X$ be a compact K\"ahler manifold, and let \holom{f}{X}{X} be an endomorphism of degree $d>1$.
Suppose that there exists a fibration $\holom{\varphi}{X}{Y}$ onto a normal variety $Y$ such that 
\begin{itemize}
\item there exists an automorphism $g$ such that $g \circ \varphi=\varphi \circ f$,
\item the general $\varphi$-fibre is isomorphic to $\PP^r$, 
\item the $\varphi$-singular locus $\Delta$ has codimension at least two in every point, and 
\item any finite \'etale cover $Y'_0 \rightarrow Y \setminus \Delta$ extends (maybe after a further finite \'etale cover) to a finite map $Y' \rightarrow Y$.
\end{itemize}
Then there exists a finite map $Y' \rightarrow Y$
such that $X \times_Y Y'$ is bimeromorphic to $Y' \times \PP^r$. 
In particular we have $a(X)=a(Y)+r$.
\end{proposition}

\begin{remark}
If $Y$ is smooth, the last condition in the proposition is automatically satisfied. In fact we have
$$
\pi_1(Y \setminus \Delta) \simeq \pi_1(Y),
$$ 
so we can even extend by an \'etale map.
\end{remark}

\begin{proof}
Since $g$ is an automorphism, the restriction of $f$ to a general fibre induces an endomorphism of degree $d>1$
and up to replacing $f$ by some $f^k$, we can suppose without loss of generality that $d>r+1$. 
Then by \cite[Prop.1.1]{Ame03} the space of endomorphism of $\PP^r$ of degree $d$ has an affine geometric quotient
$R^m(\PP^r, \PP^r)/PGL(r+1) \subset \C^N$. 
Thus the fibration $\varphi$ induces a holomorphic map $(Y \setminus \Delta) \rightarrow R^m(\PP^r, \PP^r)/PGL(r+1) \subset \C^N$. 
Since $Y$ is normal, the map extends to a holomorphical map $Y \rightarrow \C^N$ by Hartog's theorem. 
Since $Y$ is compact, this map is constant. 
Arguing as in Amerik's proof of \cite[page 22, line7ff]{Ame03}, we see that there exists an \'etale cover 
$Y'_0 \rightarrow Y \setminus \Delta$ such that the fibre product $X \times_{Y \setminus \Delta} Y'_0$ is isomorphic to $Y'_0 \times \PP^r$.
By the last condition, we know that (up to replacing  $Y'_0$ by some higher \'etale cover), the \'etale cover extends to a finite map $Y' \rightarrow Y$. 
By construction we then have a holomorphic map
$$
Y'_0 \times \PP^r \hookrightarrow X \times_Y Y',
$$
and since $(Y' \times \PP^r) \setminus (Y'_0 \times \PP^r)$ has codimension at least two, we can apply \cite[Ch.P.,Thm.10]{Gun90II} to get a bimeromorphic map
$$
Y' \times \PP^r \dashrightarrow X \times_Y Y'.
$$
\end{proof}

\subsection{Auxiliary results on compact K\"ahler surfaces}

Recall that by the classification of surfaces a compact complex surface of algebraic dimension zero 
that is in the Fujiki class is bimeromorphic to a torus or a K3 surface.
The following technical lemma is well-known to experts, but for
the convenience of the reader we include it and its (easy) proof.

\begin{lemma} \label{lemmaK3torus}
Let $S$ be a normal complex compact surface of algebraic dimension zero that is in the Fujiki class. 
Then there exists a bimeromorphic map 
$S \rightarrow S_{min}$ onto a normal surface $S_{min}$ that does not contain any curves.

If $S$ is bimeromorphic to a torus, then $S_{min}$ is a torus. 
If $S$ is bimeromorphic to a K3 surface, then  $S_{min}$ 
has at most rational double points.

Moreover if $D$ is an effective, non-trivial Cartier divisor on $S$ then $D^2<0$.
\end{lemma}

\begin{remark*}
If $S$ is bimeromorphic to a K3 surface, we will call $S_{min}$ a singular K3 surface. 
Although this has a completely different meaning in the theory of lattices of K3 surfaces, we hope
that no confusion will arise.
\end{remark*}

\begin{proof} 
Suppose first that $S$ is smooth, and 
let $S \rightarrow S'$ be the minimal model of $S$. If $S'$ is a torus,
it has no curves since the algebraic dimension is zero. 
Thus in this case we can just set $S_{min}=S'$.

If $S'$ is a K3 surface, we proceed as follows:
let $D$ be an effective divisor on $S'$. 
Since the algebraic dimension is zero, we have $h^0(S', D)=1$ and $h^2(S', D)=0$ by Serre duality, so the Riemann-Roch formula yields
$$
- h^1(S', D) = \frac{1}{2} D^2 + 1.
$$
In particular we have $D^2 \leq 2$. 
On the other hand by the adjunction formula $2 p_a(D)-2=D^2$, so the non-negativity of the arithmetic genus $p_a(D)$ implies that $D^2=-2$.
If we apply this to effective reduced divisors $D$ with one, two and three irreducible components, we see that every irreducible curve is smooth
and isomorphic to $\PP^1$, two curves are disjoint or meet transversally in exactly one point and three curves never meet in the same point.
Thus by \cite[Lemma 2.12]{BHPV04} the configuration of curves are of type A-D-E. 
Since there are only finitely many divisors on $S$ \cite[Thm.]{FF79},
we can contract all the divisors by Grauert's criterion and obtain a normal surface $S_{min}$ with at most rational double points.

If $S$ is singular we apply the first step to some desingularisation: thus there exists a meromorphic map \merom{\mu}{S}{S_{min}}.
Let $\Gamma \subset S \times S_{min}$ be the graph of the map, then the projection $\Gamma \rightarrow S$ is bimeromorphic 
and has connected fibres by Zariski's lemma. Yet any positive-dimensional fibre would be a curve in $S_{min}$, so $\Gamma \simeq S$
and $\mu$ extends to an isomorphism.

Since any effective divisor on $S$ is contracted by $\mu$ to a point, the last statement is immediate from the easy direction of Grauert's criterion.  
\end{proof}

\begin{corollary} \label{corollaryK3}
In the situation of the Lemma \ref{lemmaK3torus}, set $S_0 := S_{min} \setminus \{ p_1, \ldots, p_r \}$ and let 
$\mu_0: S'_0 \rightarrow S_0$ be an \'etale morphism. Then $\mu_0$ extends to a finite map $\mu: S' \rightarrow S_{min}$.
\end{corollary}

\begin{proof}
It is sufficient to show that we can extend $\mu_0$ locally. Yet by the lemma, the surface $S_{min}$ has at most rational double points and these
are quotient singularities of the form $\D^2/G$ where $G$ is a finite group. Thus if we take an \'etale cover of $\D^2 \setminus (0,0)/G$,
we can lift it to the universal cover $\D^2 \setminus (0,0)$ and extend by the inclusion $\D^2 \setminus (0,0) \rightarrow \D^2$.
\end{proof}

Since a smooth K3 surface is simply connected, it does not admit a (necessarily \'etale) 
endomorphism of degree at least two. This is no longer true for the singular K3 surface $S_{min}$.

\begin{example}
Let $A$ be a two-dimensional torus of algebraic dimension zero and let $S$ be the corresponding Kummer surface.
It is not hard to see that the Kummer quotient $A/\Z_2$ is the surface $S_{min}$.
The multiplication by $n \in \N$ gives an endomorphism of degree $n^4$ of the torus $A$ which descends to an endomorphism of
$A/\Z_2$ of degree at least two.
\end{example}

The following proposition shows that the example is essentially everything that can happen.

\begin{proposition} \label{propositionkummer}
Let $S$ be a singular K3 surface of algebraic dimension zero,
and let \holom{f}{S}{S} be an endomorphism of degree $d>1$.
Then there exists a Galois covering \holom{\nu}{A}{S} by a torus that is \'etale in codimension one
such that $f$ lifts to an endomorphism \holom{f_A}{A}{A} of degree $d$.
\end{proposition}

The proof is essentially a reproduction of the arguments used in \cite{Nak08} for the algebraic case.
For the reader's convenience, we give the basic ideas:

\begin{proof}  By Lemma \ref{lemmaK3torus}, the surface $S$ has only rational double points,
so it is Gorenstein and has only isolated quotient singularities. Since $S$ contains no curves, the endomorphism $f$ is necessarily 
\'etale in codimension one. Thus by \cite[Lemma 3.3.2]{Nak08} there exists a finite Galois covering $\holom{\nu}{A}{S}$ that is \'etale in
codimension one such that $A$ is smooth and $e(A)=0$. Moreover $\kappa(A) \geq 0$ and $A$ has algebraic dimension zero, so it is a torus
by the classification of compact complex surfaces. Replacing the covering $\nu$ by a suitable one, we may suppose that $\deg \nu \leq \deg \nu'$
for every Galois covering $\holom{\nu'}{A'}{S}$ by a torus $A'$ that is \'etale in codimension one. 
Arguing as in \cite[Lemma 2.6]{NZ07}, one sees that such a $\nu$ is unique up to an isomorphism over $S$.
Let $W$ be the normalisation of an irreducible component of the fibre product $A \times_S S$ such that the natural morphisms
$f': W \rightarrow A$ and $\nu': W \rightarrow S$ are surjective. Then we have a commutative diagram
$$
 \xymatrix{ 
W \ar[d]^{f'} \ar[r]^{\nu'} & S  \ar[d]^{f}  \\
A  \ar[r]^{\nu} & S}
$$
and $f'$ is \'etale in codimension one. The variety $A$ being smooth the map $f'$ is \'etale, so $W$ is a torus.
By construction, we have $\deg \nu' \leq \deg \nu$ so the minimality of $\nu$ implies that there exists 
an isomorphism $\psi: A \rightarrow W$ such that $\nu' \circ \psi=\nu$. The morphism $f_A:=f' \circ \psi$
has the stated properties.
\end{proof}

We would like to thank C. Favre for suggesting to us the proof of the next statement.

\begin{proposition} \label{propositionfavre}
Let $S$ be a compact K\"ahler surface that admits a relatively minimal elliptic fibration \holom{\varphi}{S}{C \simeq \PP^1}.
Suppose that $S$ admits a meromorphic endomorphism \merom{f}{S}{S} of degree $d>1$
such that there exists an endomorphism \holom{g}{C}{C} of degree at least two such that $g \circ \varphi=\varphi \circ f$. 
Then $S$ is algebraic.
\end{proposition}

\begin{proof} 
We argue by contradiction and suppose that $S$ has algebraic dimension one.
Since $g$ has degree at least two, there are infinitely many $\varphi$-fibres that are isomorphic elliptic curves.
Since the $j$-invariant yields a holomorphic map \holom{j}{C}{\PP^1}, it must be constant. Thus all the smooth
fibres are isomorphic and a look at the local behaviour of the $j$-invariant near the singular fibres \cite[V.10,Table 6]{BHPV04}
shows that the singular fibres are multiple elliptic curves. Since $S$ contains no curve that maps onto $C$, this implies that 
$S$ contains no rational curves. Thus the meromorphic endomorphism extends to a holomorphic endomorphism 
\holom{f}{S}{S} of degree $d>1$. Since $g$ is an endomorphism of $\PP^1$ it is ramified, so $f$ is ramified. Thus $S$ is uniruled by \cite[Thm.4.1]{AKP08}, a contradiction. 
\end{proof}

We will need the following generalisation of  \cite[Thm.4.1]{AKP08} for singular surfaces.

\begin{lemma} \label{lemmaramificationuniruled}
Let $X$ be a compact K\"ahler threefold, and let \holom{f}{X}{X} be an endomorphism of degree $d>1$.
Let $D \subset X$ be an irreducible divisor 
such that $\fibre{f}{D}=D$ and such that $\holom{f_D}{D}{D}$ has degree at least two. 
If $f_D$ is ramified, the surface $D$ is uniruled.
\end{lemma}

\begin{proof}
We claim that the canonical divisor $K_D$ is not pseudoeffective.
Assuming this for the time being, let us show how to conclude: 
let $\holom{\nu}{\tilde{D}}{D}$ be the normalisation, then 
$$
K_{\tilde{D}} = \nu^* K_D - N
$$ 
where $N$ is an effective divisor.
Let $\holom{\pi}{D'}{\tilde{D}}$ be the minimal resolution, then 
$$
K_{D'} = \pi^* K_{\tilde{D}} - N'
$$ 
where $N'$ is an effective divisor.
Thus 
$$
K_{D'} = \pi^* \nu^* K_D - N' - \pi^* N
$$ 
is not pseudoeffective, since $K_D$ is not pseudoeffective. Since $D'$ is a compact K\"ahler surface,
this implies that $D'$ is uniruled.

{\it Proof of the claim:} We start by establishing a ramification formula on $D$. Since $\fibre{f}{D}=D$, we have $f^* D = m D$ where $m \in N$ is the order of
ramification along $D$. Thus the ramification divisor $R$ of $f$ is of the form
$$
R = (m-1) D + R'
$$
where $R'$ does not contain $D$. Since $f_D$ is ramified, the restriction
$R'_D := R' \cap D$ is an effective, non-trivial Cartier divisor. By the ramification formula on $X$ we have
$K_X = f^* K_X + R$, so the adjunction formula $K_D=K_X|_D+D_D$ implies
$$
K_D - D_D = f_D^* (K_D-D) + R_D = f_D^* K_D - m D_D + (m-1) D_D + R'_D, 
$$
so
$$
K_D = f_D^* K_D +R'_D.
$$
We proceed now as in the proof of \cite[Thm.4.1]{AKP08}: let $f_m$ be the $m$-th iterate of $f$, then the ramification formula reads
$$
K_D = f_m^* K_D + f_{m-1}^* R'_D + \ldots + f^* R'_D + R'_D.
$$
Let $\omega_D$ be the restriction of a K\"ahler form $\omega$ to $D$.
There exists a $c > 0$ such that for every pseudoeffective line bundle $L$ on $D$, we have $L \cdot \omega>c$. If $K_D$ is pseudoeffective,
then $f_m^* K_D \cdot \omega \geq 0$, so
$$ 
K_D \cdot \omega \geq m c
$$
for arbitrary $m$ which is impossible. 
\end{proof}

\begin{lemma}
Let $S$ be a complex Gorenstein surface, i.e.\ the canonical sheaf exists and is locally free. Let
$\holom{f}{S}{S}$ be an endomorphism such that $K_S \simeq f^* K_S$. Let
$\holom{\nu}{\tilde{S}}{S}$ be the normalisation, and denote by $N$ the effective Weil divisor on $\tilde{S}$ such that $K_{\tilde{S}} = \nu^* K_S - N$.
There exists an endomorphism $\holom{\tilde{f}}{\tilde{S}}{\tilde{S}}$ such that $\nu \circ \tilde{f}=f \circ \nu$.
Furthermore we have $\fibre{\tilde{f}}{N} = N$ and $\tilde{f}$ is \'etale in codimension one in the complement of $N$.
\end{lemma}

\begin{proof}
The existence of $\tilde{f}$ is immediate by the universal property of the normalisation.

For the second statement, note that we have an equality of Weil divisors
$$
\tilde{f}^* (K_{\tilde{S}}+N) = \tilde{f}^* \nu^* K_S  = \nu^* f_S^* K_S
= K_{\tilde{S}} + N.
$$ 
Since by the ramification formula for normal surfaces we also have an equality of Weil divisors
$$
K_{\tilde{s}} = \tilde{f}^* K_{\tilde{S}} + R,
$$
where $R$ is the ramification divisor, we get an equality of Weil divisors
$$
R = f^* N - N,
$$
so $N$ is a completely invariant divisor of $f_{\tilde{D}}$ and  $f_{\tilde{D}}$ is \'etale outside $N$. 
\end{proof}

\section{Compact K\"ahler threefolds} \label{sectionthreefolds}

\subsection{Algebraic connectedness}

\begin{definition}
Let $X$ be a compact K\"ahler manifold. $X$ is called algebraically 
connected if there exists a family of curves $(C_t)_{t \in T}$ such that $C_t$ 
is irreducible for general $t \in T$ and 
such that two very general points can be joined by a chain of $C_t$'s.
\end{definition}
 
The following theorem of Campana illustrates the importance of this notion.

\begin{theorem} \cite{Cam81} \label{theoremcampana}
An algebraically connected compact K\"ahler manifold is projective.
\end{theorem}
 
An immediate consequence of this theorem is that if 
$X$ is a compact K\"ahler threefold of algebraic dimension two,
then the algebraic reduction of $X$ is almost holomorphic.
If the algebraic dimension is one, this no longer holds in general, but only in special cases:

\begin{corollary} \label{corollarygeneralfibrealgebraic}
Let $X$ be a compact K\"ahler non-algebraic threefold. 
Let \merom{\varphi}{X}{C} be a meromorphic map onto a curve such that
the general fibre is an algebraic surface.
Then $\varphi$ is holomorphic.
\end{corollary} 

Indeed the fibres of $\varphi$ are algebraically connected and cover $X$. Thus they can't meet since otherwise $X$
is algebraically connected. This shows that $\varphi$ is almost holomorphic, but an almost holomorphic map onto a curve
is holomorphic. 

The existence of an endomorphism of degree at least two allows us to assure the holomorphicity of the algebraic reduction in another situation.

\begin{proposition} \label{propositionexistencefibration}
Let $X$ be a compact non-algebraic K\"ahler threefold of algebraic dimension one, and let $\holom{f}{X}{X}$ be an endomorphism of degree $d>1$.
Suppose that the general fibre of the algebraic reduction is a compact K\"ahler surface of algebraic dimension zero.
Then the algebraic reduction of $X$ is holomorphic.   
\end{proposition}

\begin{proof}
Let $\holom{\mu}{X'}{X}$ be a bimeromorphic map such that $\holom{\varphi}{X'}{C}$ is a holomorphic model of the algebraic reduction. 
The endomorphism $f$ induces a meromorphic endomorphism $\merom{f'}{X'}{X'}$. 
By Proposition \ref{propositionalgebraicreduction} there exists an endomorphism \holom{g}{C}{C}
such that $g \circ \varphi=\varphi \circ f'$.
Since the general $\varphi$-fibre $F$ is irreducible, we have $\fibre{f'}{F}=F_1 \cup \ldots F_\delta$, where $\delta$ is the degree of $g$.
For simplicity denote by $F$ also the image of a general fibre under the birational map $\mu$.
Since $\fibre{f}{F}$ has $\delta$ irreducible components and $F$ does not lie in the branch locus of $f$,
we see that $f^* [F] \equiv_{num} \delta [F]$.
Thus Lemma \ref{lemmaintersection} shows that $F^3=0$ or $d=\delta^3$.

{\it 1st case. Suppose that $F^3=0$.} We want to prove that if $F'$ is another general member of the family, then $F \cap F' = \emptyset$.
Indeed if this is not the case, there exists a nontrivial effective Cartier divisor $F' \cap F$ on $F$ such that $(F'_F)^2=0$, a contradiction to  
Lemma \ref{lemmaK3torus}.
Since $F$ and $F'$ are irreducible and homologous, this implies that $[F]^2=0$. Thus a general $F$ does not meet any member of the family.
Hence the algebraic reduction of $X$ is almost holomorphic onto a curve, so it is holomorphic.

{\it 2nd case. Suppose that $F^3 \neq 0$.} In this case we have $\delta=\sqrt[3]{d}$. 
By \cite[Cor.7.6.(ii)]{CP00}, the general $\varphi$-fibres are isomorphic. Thus we can consider
the induced map $\holom{f|_{F_1}}{F_1}{F}$ as an endomorphism $\holom{f_F}{F}{F}$ of degree $\delta^2>1$. 
Since $F^3 \neq 0$, the intersection $F \cap F'$ with another general member is not trivial. The surface $F$ has only finitely many divisors \cite[Thm.]{FF79}, 
so we may suppose that $\fibre{f_F}{F \cap F'}=F \cap F'$. 
Yet $F$ is not uniruled, so the endomorphism $f_F$ is not ramified by Lemma \ref{lemmaramificationuniruled}.  
Another application of Lemma \ref{lemmaintersection} shows that $[F \cap F']^2=0$. Now $[F \cap F']^2=F \cdot F' \cdot F' = F^3$ gives a contradiction.
\end{proof}

\subsection{Mori program for compact K\"ahler threefolds} \label{subsectionmori}

Let $X$ be a compact K\"ahler manifold. A {\it contraction} is defined to be a surjective map
with connected fibres $\holom{\varphi}{X}{X'}$ onto a normal complex variety
such that $-K_X$ is relatively ample and $b_2(X)=b_2(X')+1$.
In general $X'$ might not be a K\"ahler variety, in particular $\varphi$ is not necessarily
a contraction of an extremal ray in the cone ${\overline {NE}}(X).$ 
The existence and structure of threefold contractions is assured by following statement.

\begin{theorem} \cite[Thm.2]{Pet01}, \cite[Main Thm.]{Pet98} \label{theoremcontraction}
Let $X$ be a  smooth compact K\"ahler threefold  with $K_X$ not nef.
Then $X$ carries a contraction $\holom{\varphi}{X}{X'}$ unless (possibly) $X$ is simple with $\kappa (X) = - \infty.$
The contraction is of one of the following types.
\begin{enumerate}
\item $\varphi $ is a $\PP_1$- bundle or a conic bundle over
a smooth non-algebraic surface,
\item $\varphi$ is bimeromorphic contracting an irreducible divisor
$E$ to a point, and $E$ together with its normal bundle $N_{E/X}$ is one of the
following $$(\PP_2,\sO_{\PP^2}(-1)), (\PP_2,\sO_{\PP^2}(-2)), (\PP_1 \times \PP_1, \sO_{\PP^1 \times \PP^1}(-1,-1)), (Q_0,\sO_{Q_0}(-1)),$$
where $Q_0$ is the quadric cone, 
\item $X'$ is smooth and $\varphi$ is the blow-up of $X'$ along a
smooth curve.
\end{enumerate}
The variety $X'$ is (a possibly singular) K\"ahler space in all cases except possibly 3.). 
Moreover in all cases but possibly 3.), the morphism $\varphi$ is the contraction
of an extremal ray in the cone ${\overline {NE}}(X).$ \\
In case 3.), the variety $X'$ is K\"ahler if and only if the ray of a fiber $l$ of $\varphi$ is extremal in the dual
K\"ahler cone $\overline{NA}(X).$
\end{theorem}  

The main defect of the preceding statement (and the main open problem in the Mori of compact K\"ahler threefolds) is that in general
it is not clear whether $X'$ is K\"ahler. We  give an affirmative answer in a number of situations. First we need some preparation.

\begin{lemma}  \label{currentlemma} 
Let $\psi: X \to X'$ be the blow-up of a smooth curve $B \subset X'$ in the compact manifold $X'.$ Assume
$X$ is K\"ahler. Then $X'$ is K\"ahler provided the following condition (*) is satisfied. \\
If $T$ is a positive closed current of bidimension $(1,1)$ such that $T = dS$ with some current $S,$ then $T = 0.$ 
\end{lemma}

\begin{proof} Let $b \in B$ and $l = \psi^{-1}(b). $ By Theorem \ref{theoremcontraction} we need to show that $[l]$ is extremal
in the dual K\"ahler cone $\overline{NA}(X).$ So suppose that 
$$ [l] = a_1 + a_2 $$
with $a_i \in \overline{NA}(X) \setminus 0.$ 
Now represent $a_i$ by a positive closed current $T_i$ (see e.g. \cite[Prop.1.8]{OP04}, which is a consequence of the work
of Demailly-Paun \cite{DP04}). So 
$$ l \sim T_1 + T_2, $$
and therefore
$$ \psi_*(T_1) + \psi(T_2) \sim 0.$$
By our assumption we obtain $\psi_*(T_1) + \psi_*(T_2) = 0$, so that $\psi_*(T_1) = \psi_*(T_2) = 0.$ 
Therefore the $a_i$ are proportional to $[l]$ and $[l]$ is extremal.
\end{proof} 

In the following corollary, $\chi_B$ denotes the characteristic function of $B$ and $T_B$ is the current given by
integration over $B.$ 

\begin{corollary} \label{currentcorollary} 
Let $\psi: X \to X'$ be the blow-up of a smooth curve $B \subset X'$ in the compact manifold $X'$, and denote by $E$ the exceptional divisor. Assume
$X$ is K\"ahler. Then $X'$ is K\"ahler if one of the following assertions holds.
\begin{enumerate}
\item There is no positive closed current $\lambda T_B + T'$ with $\lambda > 0$ and $\chi_E T' = 0$  which is a boundary $dS.$ 
\item A multiple of $B$ moves in a positive-dimensional family.
\item The normal bundle $N_{B/X'}$ is not negative. 
\end{enumerate}
\end{corollary} 

\begin{proof} (1) By \cite{Siu74}, see also \cite{Dem01} we have a decomposition 
$$ T = \lambda T_B + T'$$
with $T' = \chi_{X \setminus B}T$ and some $\lambda \geq 0.$ In particular $T'$ is again closed. 
Since we assume that (*) already holds for all $T$ for which 
$\lambda > 0,$ we may assume that $\lambda = 0$ and need to verify (*) for those $T.$ 
By the proof of \ref{currentlemma} we need to check that only for the current $\psi_*(T_1 + T_2). $
But for those currents $\lambda = 0$ can actually never happen: if 
$$ \chi_B \psi_*(T_1 + T_2) = 0,$$
then 
$$ \chi_E (T_1 + T_2) = 0,$$ 
where $E = \psi^{-1}(B).$ 
Thus $T_l \sim T$ with $\chi_E T = 0.$ 
Now $l \cdot E = -1;$ on the other hand $\chi_E T = 0$ implies via Demailly's regularization theorem that $E \cdot [T] \geq 0 $
\cite[Sublemma 7.a]{Pet98}.
\vskip .2cm \noindent
(2) By (1) it suffices to show the following. If $\lambda T_B + T' = dS $ with $T'$ a positive closed current such that 
$\chi_B T' = 0$ and $\lambda \geq 0$, then $\lambda = 0.$ So assume $T_B + T' = dS.$ Now a multiple of $B$ moves, so we obtain 
an irreducible curve $B' \ne B$ and $a > 0$ such that $aT_{B'} + T' = dS'.$ But $\chi_{X' \setminus B}(aT_{B'} + T') = 0,$
so the arguments of (1) show that $a = 0$ and $T' = 0$. 
\vskip .2cm \noindent
(3) Suppose that the normal bundle $N_{B/X'}$ is not negative. 
If $X'$ were not K\"ahler, then by 
\ref{currentcorollary}, we find a positive closed current $T$ with $\chi_B T = 0$ such that
$$ B + T = dS. $$
We find easily a positive closed current $\tilde T$ on $X$ such that $\chi_E \tilde T = 0$ and such that $\psi_*(\tilde T) \sim T$
(sse the proof of Prop. 2.1 in \cite{OP04}). Now choose a positive closed current $B_0$ on $E$ representing the second ray of $\overline{NE}(E)$ (i.e.\ the one not
represented by the ruling line $l$). Let $R = i_*(B_0)$ where $i: E \to X$ is the inclusion. Then up to scaling $\psi_*(R) = B$ and therefore
$$ R + \tilde T \sim al $$
for some $a > 0.$ Intersecting with $E$ and observing $\tilde T \cdot E = 0 $ as in the proof of \ref{currentcorollary} we obtain $B_0 \cdot E < 0$ so that 
the normal bundle $N_{E/X}$ is negative and so does $N_{B/X'}. $ 
\end{proof} 

We can now prove that the study of compact K\"ahler threefolds with endomorphisms can be reduced to minimal models,
at least if $X$ is not uniruled.

\begin{proposition} \label{propositionmmpnotuniruled}
Let $X$ be a compact K\"ahler threefold which is not uniruled, 
and let \holom{f}{X}{X} be a (necessarily \'etale) endomorphism of degree $d>1$.
Then there exists a finite sequence
of smooth compact {\em K\"ahler} threefolds
\[
X=X_0 \stackrel{\mu=\mu_1}{\longrightarrow} X_1 \stackrel{\mu_2}{\longrightarrow} \ldots \stackrel{\mu_n}{\longrightarrow} X_n
\]
such that 
\begin{itemize}
\item the canonical bundle of  $X_n$ is nef;
\item for all $i \in \{1, \ldots, n \}$, the manifold $X_{i-1}$ is the blow-up of $X_i$ along a smooth elliptic curve $C_i \subset X_i$ such that
the normal bundle is an indecomposable rank two vector bundle of degree zero or
a direct sum of numerically trivial line bundles;
\item up to replacing $f=:f_0$ by $f^k$ for some $k \in \N$, there exist endomorphisms $\holom{f_i}{X_i}{X_i}$ of degree $d$
such that $f_i \circ \mu_i=\mu_i \circ f_{i-1}$ for all $i \in \{1, \ldots, n \}$. The curve $C_i$ is invariant under $f_i$, i.e.\ we have $\fibre{f_i}{C_i}=C_i$.
\end{itemize} 
\end{proposition}

\begin{proof} We will use Theorem \ref{theoremcontraction}: note first that 
the case where $X$ is simple with $\kappa (X) = - \infty$ can be easily excluded in our situations:
since $f$ is \'etale by Proposition \ref{propositionuniruled} we have $\chi(X, \sO_X)=0$ which implies that $h^0(X, K_X) \geq 1$ or $h^0(X, \Omega_X) \geq 1$. In
the first case we are obviously done and in the second case we can look at the Albanese map $X \rightarrow Alb(X)$ and conclude
that $\kappa(X) \geq 0$ by the $C_{n,m}$-conjecture (which is known for compact K\"ahler threefolds \cite[Thm.2.2, Thm.4.1]{Uen87}).

Thus there exists a $m \in \N$ such that $|m K_X|$ is not empty. 
Let $D$ be the fixed part of this linear system, then $f^* K_X \simeq K_X$ implies that $f^* D \simeq D$. 
Since all the contractions are divisorial by Theorem \ref{theoremcontraction}, their exceptional loci are contained in $D$.
Since $D$ has only finitely many components, the description of the exceptional loci implies that there are only finitely many
contractions. The endomorphism $f$ permutes the irreducible components of $D$, so up to replacing $f$ by $f^k$ we may suppose
that $\fibre{f}{E}=E$ for the support $E$ of every birational contraction.

Choose now a birational contraction \holom{\mu}{X}{X_1} with exceptional divisor $E$.
By Lemma \ref{lemmaintersection} we have $E^3=0$, so the description of the normal bundles in 
Theorem \ref{theoremcontraction} shows that $E$ is contracted onto a smooth curve $C_1$ and $X_1$ 
is a smooth complex manifold.
The restriction of $f$ to $E \simeq \PP(N^*_{C_1/X_1})$ 
yields an \'etale endomorphism \holom{f}{E}{E} of degree $d$, so \cite[Thm.1.1.]{Nak08} shows that $C_1$ is elliptic
and $N^*_{C_1/X_1}$ has the prescribed form, at least up to twist by a line bundle. Since
$$
c_1(N^*_{C_1/X_1})=\sO_{\PP(N^*_{C_1/X_1})}(-1)^2=E^3=0,
$$ 
we see that such a line bundle is numerically trivial.
This completes the description of the normal bundle of $C_1$. 
Corollary \ref{currentcorollary} now implies that $X_1$ is a K\"ahler manifold.
By \cite[Rem.3.3,Cor.3.4]{AKP08} there exists a $k_1 \in \N$ and an
an endomorphism $\holom{f_1}{X_1}{X_1}$ of degree $d$
such that $f_1 \circ \mu_1=\mu_1 \circ f$. 

Hence $X_1$ satisfies all the assumption made on $X$ and we can argue by induction on the Picard number to get the existence of the sequence satisfying the
first two properties. Replacing the $k_i$ by a sufficiently divisible $k$, we obtain the third statement.
\end{proof}

If $X$ is uniruled the situation gets more complicated. One reason is that the exceptional divisor of a birational contraction can lie
in the branch locus of $f$. If this is not the case, we can again use 
Lemma \ref{lemmaintersection}, Theorem \ref{theoremcontraction} and Corollary \ref{currentcorollary} to show the next statement.

\begin{corollary} \label{corollarycontractionuniruled}
Let $X$ be a compact K\"ahler threefold, and let \holom{f}{X}{X} be an endomorphism of degree $d>1$.
Let \holom{\psi}{X}{X'} be a contraction of birational type such that the exceptional divisor $E$ satisfies $\fibre{f}{E}=E$.
Suppose that $E$ is not contained in the branch locus of $f$ or that the ramification order along $E$ is not $\sqrt[3]{d}$.  
Then $E^3=0$ and $\psi$ contracts the divisor $E$ onto a smooth curve $C'$. Moreover the variety $X'$ is a compact K\"ahler threefold.
\end{corollary}

\begin{remark*}
If the endomorphism $f$ is \'etale one proves as above that the curve $C'$ is elliptic. If $f$ is ramified, $C'$ can a-priori be arbitrary.
\end{remark*}

\begin{proposition} \label{propositionkaehlerpackage}
Let $X$ be a compact non-algebraic K\"ahler threefold, and let \holom{f}{X}{X} be an endomorphism of degree $d>1$.
Suppose that there exists a fibration $\holom{\varphi}{X}{C}$ onto a smooth curve $C$ with algebraic general fiber and an  
automorphism \holom{g}{C}{C} such that $g \circ \varphi=\varphi \circ f$.

Suppose that $X$ admits a birational contraction \holom{\psi}{X}{X'}. Then 
\begin{itemize}
\item the exceptional locus $E$ is contained in a $\varphi$-fibre,
\item there exists a fibration \holom{\varphi'}{X'}{C} such that $\varphi=\varphi' \circ \psi$,
\item $\psi$ is the blow-up of a smooth curve $B$ in $X'$, 
\item there exists an endomorphism  \holom{f'}{X'}{X'} of degree $d$ such that $f' \circ \psi=\psi \circ f$, and
\item $X'$ is a compact K\"ahler threefold unless possibly when the curve $B$ has negative normal bundle.
\end{itemize}
 If $\varphi$ is locally projective, there exists a possibly different contraction $\psi': X \to X^*  $ with locally projective factorization $X^* \to C$.
\end{proposition}

\begin{proof}
By the classification in Theorem \ref{theoremcontraction}, the exceptional loci of birational contractions are uniruled surfaces and therefore algebraic.
Therefore $E$ is contained in a $\varphi$-fibre, since otherwise $X$ is algebraically connected, hence algebraic by Theorem \ref{theoremcampana}.
The existence of the fibration $\varphi'$ is now an immediate consequence of the rigidity lemma.

Since $E$ is contained in a fibre, it is not contained in the branch locus by Lemma \ref{lemmafibrationcurve}. Therefore Corollary 
\ref{corollarycontractionuniruled} shows that $E$ is contracted onto a smooth curve and $X'$ is smooth. 
By what precedes we know that the exceptional loci of birational contractions on $X$ are irreducible components of $\varphi$-fibres. Since there are only
finitely many reducible fibres and these have finitely many components, we see that there are only finitely many birational contractions. Thus up to replacing $f$
by $f^k$ we have $f_* [\Gamma]=\lambda [\Gamma]$, where $\Gamma$ is a fibre of $E \rightarrow \psi(E)$. We conclude as in \cite[Cor.3.4]{AKP08},
the K\"ahler property follows from Corollary \ref{currentcorollary}.

Thus we are left with the last assertion. Let                                             
$$\{x_1, \ldots, x_s\} \subset C$$ 
be the singular locus of $\varphi$.
Choose small disjoint open neighborhoods $U_i$ of $x_i$ such that $\varphi$ is projective over $U_i.$ 
Let $X_i = \varphi^{-1}(U_i)$ and $\varphi_i = \varphi \vert X_i.$ 
We already know that there exists some $i$ such that $K_X $ is not $\varphi_i-$nef. Hence by \cite{Nak87} there exists a relative contraction                          
$\mu: X_i \to X'_i$ such that the induced map $X'_i \to U_i$ is projective. Now $\mu $ might be birational or not. We always choose $\mu$ birational unless it is simply
not possible, i.e.\ for all choice of $i$, the map $\mu$ is a fibration.                                                                                                 
In case $\mu$ is birational we patch things and obtain $X^* $ with a 
locally projective bimeromorphic map $\varphi': X^* \to C$. The map $\mu$ extends to $\psi: X \to X^*.$ \\
In case $\mu$ is a fibration, we obtain by deformation of the extremal rational curves and by our assumption a global relative contraction which is a $\PN_1-$bundle or
a conic bundle $X \to X^*$ with factorization $X^* \to C.$                 
\end{proof}

Finally let us recall that abundance holds for most minimal K\"ahler threefolds.

\begin{theorem} \cite[Thm.1]{Pet01} \label{theoremabundance}
Let $X$ be a normal compact K\"ahler threefold (${\mathbb Q}-$factorial with
at most terminal singularities) such that $K_X$ is nef. Assume that $X$ is not both simple and non-Kummer. Then
$K_X$ is semi-ample. 
\end{theorem}

\section{Torus fibrations} \label{sectiontorusfibrations}

\begin{definition}
A torus fibration is a fibration \holom{\varphi}{X}{Y} such that the general fibre is isomorphic to a complex torus.
A torus bundle is a smooth torus fibration that is locally trivial.
\end{definition}

If the total space of a torus fibration is not projective, the fibration in general does not admit a multisection, i.e.\ there is no subvariety $Z \subset X$
such that \holom{\varphi|_Z}{Z}{Y} is surjective and generically finite. The main technical statement of this section (Lemma \ref{lemmamultisection}) shows that if $X$
admits an endomorphism commuting with $\varphi$, then there exists a natural meromorphic factorisation of $\varphi$
which admits a multisection.

\begin{lemma} \label{lemmamultisection}
Let $X$ be a compact normal variety in the Fujiki class that admits a torus fibration \holom{\varphi}{X}{Y}.
Suppose furthermore that there exists an endomorphism \holom{f}{X}{X} of degree $d>1$
and an automorphism \holom{g}{Y}{Y} of finite order such that $g \circ \varphi=\varphi \circ f$.
Then (up to replacing $f$ by some power) there exists a compact normal variety $Z$ that is in the Fujiki class 
and admits a torus fibration \holom{\psi}{Z}{Y} with a multisection and which satisfies the following properties:
\begin{itemize}
\item there exists an almost holomorphic fibration \merom{\tau}{X}{Z} such that $\varphi=\psi \circ \tau$.
\item there exists an endomorphism \holom{\overline{f}}{Z}{Z} of degree $d$ that commutes with $\psi$ and such that
$\overline{f} \circ \tau=\tau \circ f$. 
\end{itemize}
If a very general fibre of $\varphi$ is a simple torus, then $\varphi$ admits a multisection.
\end{lemma}

The proof of this lemma is based on the following easy observation.
The statement generalises \cite[Lemma 2.22]{FN07} whose
strategy of proof we follow.

\begin{proposition} \label{propositiontorusfixedpoints}
Let $A$ be a complex torus, and let $\holom{f}{A}{A}$ be an endomorphism of degree $d>1$. 
Then there exists a (maybe trivial) subtorus $T \subsetneq A$
and an endomorphism
$\holom{\overline{f}}{A/T}{A/T}$ of degree $d$ such that the set of fixed points of $\overline{f}$ is non-empty and finite.

In particular if $A$ is simple,  the set of fixed points of $f$ is non-empty and finite.
\end{proposition}     

\begin{proof}[Proof of Proposition \ref{propositiontorusfixedpoints}]
We will argue by induction on the dimension, the case of dimension one is included in \cite[Lemma 2.22]{FN07}.
Choose a point $0 \in A$ so that the torus $A$ has a group structure. 
The map 
$$
\holom{h}{A}{A}, x \mapsto f(x)-f(0)-x
$$
is a morphism of groups and not zero, since $f$ has degree at least two. Let $T_0$ be the connected component of the kernel of $h$.
We make a case distinction.

{\em 1st case. $T_0$ is trivial} In this case $h$ is surjective and has finite kernel. Since
$$
\{ x \in A \ | \ f(x)=x \} = \{ x \in A \ | \ h(x)=f(0) \} 
$$ 
the statement follows.

{\em 2nd case. $T_0$ has positive dimension.} 
Since $h(x)=0$ for all $x \in T_0$, we have $f(x)=x+f(0)$ for all $x \in T_0$.
It is thus clear that there exists an endomorphism
$\holom{\overline{f}}{A/T_0}{A/T_0}$ of degree $d$ such that $\overline{f} \circ q=q \circ f$, where $\holom{q}{A}{A/T_0}$
is the quotient map. Apply the induction hypothesis to $A/T_0$.
\end{proof}

\begin{example}
At first glance the proof of Proposition \ref{propositiontorusfixedpoints} may suggest that the restriction
of the endomorphism $f$ to $T$ is a translation. The endomorphism of degree $n$
$$
f: E \times E \times E \rightarrow E \times E \times E, \ (x_1, x_2, x_3) \ \mapsto (x_1+x_2+x_3, x_2+x_3, n x_3) 
$$
shows that this is not true, since we will have $T=E \times E \times \{ 0 \}$. The restriction of $f$ to $T$ is rather a ``tower''
of translations. 
\end{example}

\begin{proof}[Proof of Lemma \ref{lemmamultisection}]
The automorphism $\holom{g}{Y}{Y}$ is assumed to be of finite order, so up 
to replacing $f$ by some multiple we can suppose that $g=Id_Y$.
The statement claims that there exists a commutative diagram 
$$
 \xymatrix{ 
X \ar @{.>}[rd]^\tau \ar[ddr]_{\varphi} \ar [rrr]_{f'} & &  &  X  \ar @{.>}[ld]_\tau \ar[ddl]^{\varphi}
\\
& Z \ar[d]^{\psi} \ar[r]^{\overline{f}} & Z  \ar[d]^{\psi} &  
\\
& Y  \ar[r]^{Id_Y} & Y & 
}
$$
such that $\overline{f}$ has degree $d$ and $Z \rightarrow Y$ has a multisection.

We denote by $\holom{f_y}{X_y}{X_y}$ the restriction of $f$ to a general fibre $X_y$: this is an endomorphism of degree $d$.
If the fix point set of $f_y$ is finite, the fix point set of $f$ is a multisection of $\varphi$, thus the statement is trivially true.
Note that if $X_y$ is simple, we are always in this case, which proves the last part of the statement. 

Suppose now that the fix point set of $f_y$ is not finite.
Then there exists by Proposition \ref{propositiontorusfixedpoints} a torus $T_y \subset X_y$ and an 
endomorphism \holom{\overline{f}_y}{X_y/T_y}{X_y/T_y} of degree $d$ that commutes with the projection
$X_y \rightarrow X_y/T_y$ and has non-empty finite fix point set.
The countability of the number of irreducible components of the relative cycle space $\chow{X/C}$ implies
that if we choose $y$ very general,
the torus $T_y$ deforms with $X_y$. 
Let $Z \rightarrow Y$ be the normalisation of the component of $\chow{X/Y}$ 
whose very general points corresponds to 
the translates of $T_y$ for $y \in Y$ very general.
The almost holomorphic map \merom{\tau}{X}{Z} such that $\varphi=\psi \circ \tau$ is given fibrewise 
by $X_y \rightarrow X_y/T_y$.

 The (necessarily \'etale) endomorphism $f$ acts via the push-forward of cycles on $\chow{X/Y}$.
Noting that the general fibre of $Z \rightarrow Y$ is the quotient torus $X_y/T_y$ and
the restriction of $f_*$ to $X_y/T_y$ identifies to $\overline{f}_y$, we see that $f_*$ maps the component 
corresponding to $Z$ onto itself. Thus we obtain a holomorphic
endomorphism \holom{\overline{f}}{Z}{Z}
of degree $d$ that commutes with $\psi$. The restriction of  $\overline{f}$ to a very general fibre $Z_y$
has non-empty finite fix point set, so the fix point set of $\overline{f}$ is a multisection of $\psi$.
By construction of $\overline{f}$, it is clear that $\overline{f} \circ \tau=\tau \circ f$.
\end{proof}

Theorem \ref{theoremkappanminusone} is now an immediate application:

\begin{proof}[Proof of Theorem \ref{theoremkappanminusone}] \label{prooftheoremkappanminusone}
We argue by contradiction and suppose that $X$ is not projective. Since $a(X) \geq \kappa(X)=n-1$, we see that the algebraic dimension is $n-1$.
By Campana's theorem \ref{theoremcampana} this implies that the algebraic reduction $\merom{\varphi}{X}{Y}$ is almost holomorphic. Since $\kappa(X)=n-1$
it is bimeromorphicly equivalent to the Iitaka fibration.  
By \cite[Thm. A]{NZ07} the endomorphism $f$ induces an automorphism $\holom{g}{Y}{Y}$ of finite order.
Thus Lemma \ref{lemmamultisection} implies that $X \rightarrow Y$ admits a multisection.
Since $Y$ is projective and $X_y$ a curve, $X$ is algebraically connected.
Hence it is projective by Campana's theorem, a contradiction.
\end{proof}  

\begin{corollary} \label{corollarymultisection}
Let $X$ be a compact K\"ahler manifold that is a torus bundle \holom{\varphi}{X}{Y} with fibre $A$.
Suppose furthermore that there exists an endomorphism \holom{f}{X}{X} of degree $d>1$
and an automorphism \holom{g}{Y}{Y} of finite order such that $g \circ \varphi=\varphi \circ f$.
Then (up to replacing $f$ by some power) there exists a compact K\"ahler manifold $Z$
that is a torus bundle \holom{\psi}{Z}{Y} with fibre $A/T$ with a multisection and satisfies the following properties:
\begin{itemize}
\item $X$ is a torus bundle \holom{\tau}{X}{Z} with fibre $T$ such that $\varphi=\psi \circ \tau$.
\item There exists an endomorphism \holom{\overline{f}}{Z}{Z} of degree $d$ that commutes with $\psi$ and such that
$\overline{f} \circ \tau=\tau \circ f$. 
\item The multisection is given by $Fix(\overline{f})$ which is an \'etale cover of $Y$.
\end{itemize}
If $A$ is simple, then $\varphi$ admits an \'etale multisection.
\end{corollary}

\begin{proof}
The proof of the first two points is the same as for Lemma \ref{lemmamultisection}, with one difference: since $X \rightarrow Y$ is a bundle,
it is clear that for very general $y$ the torus $T_y \subset A$ does not depend on $y \in Y$. 
Thus the variety $Z$ can be directly defined as the quotient bundle $Z \rightarrow Y$ with fibre $A/T$.

The \'etaleness of $Fix(\overline{f}) \rightarrow Y$ is shown by copying word by word the proof of \cite[Thm.2.24]{FN07}.
\end{proof}

Let now $X$ be a complex torus that contains a proper subtorus $T \subsetneq X$. Then $X$ is naturally a torus bundle over $X/T$ with fibre $T$
and if $X$ is projective Poincar\'e's irreducibility lemma shows that it is isogenous to $X/T \times T$. 
In particular the study of endomorphisms of $X$ is reduced to $X/T \times T$. 
If $X$ is not projective the situation is more complicated, but the geometric intuition still says that if $X \rightarrow X/T$
admits an  ``interesting'' \ endomorphism, the torus $X$ should be close to being a product. 
We illustrate this philosophy in two special cases:                                                                      

\begin{proposition} \label{propositiontorusoverelliptic}
Let $X$ be a torus that is a torus bundle $\holom{\varphi}{X}{Y}$ over a torus $Y$ with fibre $A$. 
Suppose that $\dim Y=1$ or $\dim Y=2, a(Y)=1$. 
Suppose that
$X$ admits an endomorphism \holom{f}{X}{X} of degree $d>1$ such that there exists an automorphism 
\holom{g}{Y}{Y} such that $g \circ \varphi=\varphi \circ f$. 
Then there exists a (maybe trivial) subtorus $T \subsetneq A$
such that $X$ is isogenous to a torus bundle over $E \times A/T$.

In particular if $A$ simple, then $X$ is isogenous to $Y \times A$.
\end{proposition}

\begin{proof} Choose a point $0 \in X$ so that the tori $X, Y$ and $A$ have a group structure. 
Up to composing $f$ with the translation $x \rightarrow x - f(0)$, 
we can suppose that $f(0)=0$. Thus the automorphism $g$ satisfies $g(0)=0$, i.e.\ is a group automorphism of $Y$.
Since $Y$ is a curve or a surface of algebraic dimension one, this group is finite (cf. \cite[Prop.3.10]{Fuj88} for the surface case), 
so $g$ is of finite order.
By Corollary \ref{corollarymultisection}, there exists a quotient $X \rightarrow X/T \rightarrow Y$
such that $f$ descends to an endomorphism $\overline{f}$ on $X/T$. 
Moreover the covering $Fix(\overline{f}) \rightarrow Y$ is \'etale,
so an irreducible component of $Fix(\overline{f})$ is a complement 
to $A/T$ in $X/T$. The statement follows by \cite[Prop.6.1]{BL99}.  
\end{proof}

\begin{proposition} \label{propositionellipticovertorus}
Let $X$ be a torus that is a torus bundle $\holom{\varphi}{X}{Y}$ over a torus  $Y$ with fibre $A$. 
Suppose that $End(A) \simeq \Z$ or $\dim A=1$.
Suppose that
$X$ admits an endomorphism \holom{f}{X}{X} of degree $d>1$ such that there exists an automorphism 
\holom{g}{Y}{Y} such that $g \circ \varphi=\varphi \circ f$. 
Then $X$ is decomposable.
\end{proposition}

\begin{proof}
We prove the statement in the case where $End(A) \simeq \Z$ (the same strategy works 
if $A$ is an elliptic curve with complex multiplication, cf. the proof of \cite[Prop.3.10]{Fuj88}).
We argue by contradiction and suppose that $X$ is indecomposable.
Since $End(A) \simeq \Z$, the restriction of $f$ to $A$ is the multiplication by an integer $n$ such that $n^{\dim A}=d$.
Thus $f-n$ is an endomorphism of $X$ whose restriction to $A$ is constant, in particular $f-n$ is not an isogeny.
Thus by \cite[Prop.7.3]{BL99}, the endomorphism $f-n$ is nilpotent. 
Hence the induced endomorphism $g-n$ on $Y$ is nilpotent, so its kernel has positive dimension. 
Thus there exists a positive dimensional subtorus $T \subset Y$ such that the restriction
of $g$ to $T$ equals the multiplication by $n$. In particular $g|_T$ is not injective, 
so $g$ is not an automorphism.  
\end{proof}

Combining Proposition \ref{propositiontorusfixedpoints} with Proposition \ref{propositionellipticovertorus},
we obtain:

\begin{corollary} \label{corollarytorusaonefixedpoints}
Let $A$ be a two-dimensional torus of algebraic dimension one, 
and let $\holom{f}{A}{A}$ be an endomorphism of degree $d>1$. 
Then $f$ has a non-empty finite set of fixed points.
\end{corollary}

\section{Non-uniruled manifolds} \label{sectionnonuniruled}

In this section we prove Theorem \ref{maintheoremA}. Using the minimal model program 
for non-uniruled threefolds admitting endomorphisms established in Proposition \ref{propositionmmpnotuniruled}
the proof  naturally splits into two parts: 
the first and most difficult part is to 
classify the minimal models admitting endomorphisms (Theorem \ref{theoremminimalmodels} below).
Based on the rather short list obtained in the first step, we then discuss the structure of the blow-ups
in Subsection \ref{subsectionnonminimal}.

\subsection{Minimal models}

\begin{theorem} \label{theoremminimalmodels}
Let $X$ be a smooth compact non-algebraic K\"ahler threefold which is not uniruled. 
Suppose that $X$ admits a (necessarily \'etale) endomorphism $f: X \to X$ of degree $d > 1.$ 
If $K_X$ is nef, then (up to \'etale cover) one of the following holds:
\begin{enumerate}
\item $\kappa(X)=0:$ then either 
\begin{enumerate}
\item $X$ is a torus or 
\item $X$ is a product $S \times E$ where $S$ is a non-algebraic K3 surface
and $E$ an elliptic curve.
\end{enumerate}
\item $\kappa(X)=a(X)=1:$ then $X$ is a product $C \times A$ where $C$ is a curve of general type and $A$ a torus
of algebraic dimension zero.
\item $\kappa(X)=1, a(X)=2:$  then either
\begin{enumerate}
\item $X$ is a product  $Y \times A$
where $Y$ is of general type and $A$ a torus of algebraic dimension one or
\item $X$ is a product $E \times S$
where $E$ is an elliptic curve and
and $S$ a non-algebraic K\"ahler surface of Kodaira dimension one.
\end{enumerate}
\end{enumerate}
\end{theorem}

\begin{proof}
By Theorem \ref{theoremkappanminusone} we have $\kappa(X) \leq 1$, and
by Theorem \ref{theoremabundance} the canonical bundle is semi-ample.

If $\kappa (X) =  0$ this implies $mK_X = \sO_X$ for some $m > 0$.  By the Beauville-Bogomolov decomposition theorem 
$X$ admits a finite \'etale cover by a torus, or a product of an elliptic curve and a K3 surface, or a Calabi-Yau manifold
of dimension three. The last case is excluded since a Calabi-Yau manifold is simply connected.

Thus we are reduced to study the cases $\kappa (X)=1$. These are dealt with in the Theorems \ref{theoremkappanminustwoA},
and \ref{theoremkappanminustwoB}.
\end{proof} 

\begin{lemma} \label{lemmarelativedimensiontwo}
Let $X$ be a non-algebraic compact K\"ahler manifold of dimension $n$,
and let \holom{f}{X}{X} be an endomorphism of degree $d>1$. 
Suppose that $X$ admits a fibration \holom{\varphi}{X}{Y} onto a projective variety $Y$ 
such that $\varphi \circ f=\varphi$ and the general fibre $X_y$ has Kodaira dimension zero.
Then the general fibre $X_y$ is (up to \'etale cover) a two-dimensional torus. 
\begin{enumerate}
\item If $a(X)=n-2$, the fibre $X_y$ has algebraic dimension zero and is isomorphic to fixed torus $A$.
\item If $a(X)=n-1$ and $a(X_y)=1$, the fibre $X_y$ is isomorphic to fixed torus $A$.
\item If $a(X)=n-1$ and $a(X_y)=2$, there exists a normal projective variety $Z$ 
that admits an elliptic fibration \holom{\psi}{Z}{Y} and satisfies the following properties:
\begin{itemize}
\item there exists an almost holomorphic fibration \merom{\tau}{X}{Z} such that $\varphi=\psi \circ \tau$.
\item there exists an endomorphism \holom{\overline{f}}{Z}{Z} of degree $d$ that commutes with $\psi$ and such that
$\overline{f} \circ \tau=\tau \circ f$. 
\end{itemize}
\end{enumerate}
\end{lemma}

\begin{proof}
The general fibre $X_y$ has Kodaira dimension zero and the restricted endomorphism \holom{f_y}{X_y}{X_y}
has degree $d$. Thus $X_y$ is covered by a torus by \cite{FN05}.
Moreover if $a(X_y) \leq 1$, the general fibres are isomorphic to a fixed torus $A$:
apply \cite[Cor.6.8]{CP00} to $X \times_Y C \rightarrow C$ where $C$ is a general complete intersection curve in $Y$.

Suppose now that $a(X_y)=1$. Let $\holom{r_y}{X_y}{E_y}$ be the algebraic reduction, then there exists an endomorphism
$\holom{\overline{f}_y}{E_y}{E_y}$ such that $\overline{f}_y \circ r_y = r_y \circ f_y$. The endomorphism $\overline{f}_y$ can't be an automorphism since
otherwise $a(X_y)=2$ by Proposition \ref{propositiontorusoverelliptic}. Thus it has degree at least two
and if $\psi: Z \dashrightarrow Y$ denotes the relative algebraic reduction of $X \rightarrow Y$, there exists a meromorphic
endomorphism \merom{\overline{f}}{Z}{Z} that commutes with $\psi$.
By Lemma \ref{lemmamultisection} this implies that $\psi$ has a multisection, hence
$a(X) \geq a(Z) =n-1$. 

Suppose now that $a(X_y)=2$.
Then the fixed point set of $f_y$ is not finite, since otherwise $f$ has a multisection and
$X$ is projective by Campana's theorem \ref{theoremcampana}. 
Thus by Lemma \ref{lemmamultisection} 
there exists a compact normal variety $Z$ that is in the Fujiki class and 
admits an elliptic fibration \holom{\psi}{Z}{Y} with a multisection and satisfies the stated properties.
Since $\psi$ has a multisection, we have $a(X) \geq a(Z) =n-1$ and $Z$ is algebraic.
\end{proof}

\begin{theorem} \label{theoremkappanminustwoA}
Let $X$ be a compact K\"ahler manifold of dimension $n$,
and let \holom{f}{X}{X} be an endomorphism of degree $d>1$. 
Suppose that 
$$
\kappa(X)=a(X)=n-2.
$$ 
Then $X$ is (up to \'etale cover) a product $Y \times A$
where $Y$ is of general type and $A$ a two-dimensional torus of algebraic dimension zero.
\end{theorem}

\begin{proof} 
By \cite[Sect. 1.4]{NZ09} there exists a ``$f$-equivariant'' resolution of the indeterminacies \holom{\mu}{X'}{X}
of the Iitaka fibration \merom{\varphi}{X}{Y} such that $f$ lifts to a {\em holomorphic} endomorphism
\holom{f'}{X'}{X'}. If we show that $X' \simeq Y' \times A$, it is {\em a-posteriori} 
clear that the Iitaka fibration of $X$ is holomorphic and hence $X \simeq Y \times A$.
Thus we can suppose without loss of generality that $X$ admits a holomorphic
fibration $\holom{\varphi}{X}{Y}$ onto a projective variety $Y$ 
such that the general fibre has Kodaira dimension zero.
By \cite[Thm. A]{NZ09} the endomorphism $f$ induces an automorphism $\holom{g}{Y}{Y}$ of finite order, so up 
to replacing $f$ by some multiple we can suppose that $g=Id_Y$.
By Lemma \ref{lemmarelativedimensiontwo} the general fibre $X_y$ is a torus of algebraic dimension one. 
Since by hypothesis $X_y$ is not algebraic, it is isomorphic to a fixed torus $A$ (ibid).

Thus by \cite[Cor. 6.6]{CP00} there exists a finite Galois cover $Y' \rightarrow Y$ such that 
$X \times_Y Y'$ is bimeromorphic over $Y'$ to a principal torus bundle $\varphi': X' \rightarrow Y'$.
The following commutative diagram and the universal property of the fibre product show that $f$ lifts to an endomorphism
$h: X \times_Y Y' \rightarrow X \times_Y Y'$ of degree $d$ that commutes with the projection on $Y'$.
$$
 \xymatrix{ 
X \times_Y Y' \ar[rd] \ar[ddd] \ar @{.>}[rrr]_h^{\exists \ \mbox{\small univ. prop. fibre prod.}} & &  &  X \times_Y Y'  \ar[ld] \ar[ddd]
\\
& X \ar[d]^{\varphi} \ar[r]^{f} & X  \ar[d]^{\varphi} &  
\\
& Y  \ar[r]^{Id_Y} & Y & 
\\
Y' \ar[ru] \ar[rrr]^{Id_{Y'}} & & & Y' \ar[lu]
}
$$
Thus there exists a meromorphic endomorphism $\merom{f'}{X'}{X'}$ of degree $d$ that commutes with $\varphi'$.
The $\varphi'$-fibre is a two-dimensional torus of algebraic dimension zero, so it contains no curves.
This immediately implies that $f'$ extends to a holomorphic endomorphism. 
Since $X_y$ is simple, Corollary \ref{corollarymultisection} shows that
the fixed point set of $f$ gives an \'etale covering
$Fix(f) \rightarrow Y'$.
Since $X' \rightarrow Y'$ is a principal bundle, we see that after \'etale base change
$X' \simeq Y' \times A$
Moreover since $Y'$ is projective any morphism from $Y$ to $A$ is constant, so
$f'=(Id_{Y'}, f_A)$ where $f_A$ is the restriction of $f$ to any $y \times A$.
 
Since $A$ contains no curves, the 
composed map $\merom{\psi:=p_A \circ \mu}{X \times_Y Y'}{A}$ extends to a holomorphic map.
Moreover $f_A \circ \psi=\psi \circ h$, so 
Proposition \ref{propositionspecialfibres} implies that $\psi$ is smooth.
The natural maps $X \times_Y Y' \rightarrow Y'$ and $\psi$ then define an isomorphism onto $Y' \times A$.
\end{proof}

\begin{theorem} \label{theoremkappanminustwoB}
Let $X$ be a compact K\"ahler threefold,
and let \holom{f}{X}{X} be an endomorphism of degree $d>1$. 
Suppose that $K_X$ is semiample and
$$
\kappa(X)=1, a(X)=2.
$$
Then (up to \'etale cover)
one of the following holds:
\begin{enumerate}
\item 
$X$ is a product $C \times A$
where $C$ is a curve of general type and $A$ a two-dimensional torus of algebraic dimension one or
\item $X$  is a product $E \times S$
where $E$ is an elliptic curve and $S$ a non-algebraic K\"ahler surface of Kodaira dimension one.
\end{enumerate}
\end{theorem}

The stategy of the proof should also work for compact K\"ahler manifolds
with $\kappa(X)=n-2, a(X)=n-1$ if one is able to show the
following statement (due to Nakayama \cite[Thm.6.2.1]{Nak08} in the surface case).

\begin{conjecture}
Let $X$ be a normal projective variety of dimension $n-1$,
and let \holom{f}{X}{X} be an endomorphism of degree $d>1$. 
Suppose that $X$ admits a (flat?) fibration \holom{\tau}{X}{Y} whose general fibre is an elliptic curve
and commutes with $f$.
Then (up to base change) we have $X$ is a product $Y \times E$ where $E$ an elliptic curve.
\end{conjecture}

\begin{proof}
Some multiple bundle of the canonical bundle induces a fibration
\holom{\varphi}{X}{C} such that $m K_X \simeq \varphi^* L$.
By \cite[Thm. A]{NZ09} the endomorphism $f$ induces an automorphism $\holom{g}{C}{C}$ of finite order, so up 
to replacing $f$ by some multiple we can suppose that $g=Id_C$.

{\em 1st case. The general $\varphi$-fibre has algebraic dimension one.}

By a statement due to Voisin (cf. \cite[Prop.3.10]{Cam06}), there exists a holomorphic two-form
on $X$ whose restriction to the general fibre $X_c$ gives a non-zero holomorphic two form.
Thus by \cite[Prop.4.2, Prop.6.7]{CP00} the fibration $\varphi$ is almost smooth 
and there exists a base change $C' \rightarrow C$ such that the normalisation of $X \times_C C'$ gives an \'etale covering of  $X$ 
and is a principal bundle over $C'$. Arguing as in the proof of Theorem \ref{theoremkappanminustwoA} we see that the endomorphism $f$
lifts, so we can suppose without loss of generality that $\varphi: X \rightarrow C$ is a principal bundle.
By Corollary \ref{corollarytorusaonefixedpoints} the induced endomorphism \holom{f}{X_c}{X_c} has a non-empty finite fixed point set.
Copying word by word the proof of \cite[Thm.2.24]{FN07}, one sees that $Fix(f) \rightarrow C$
is finite and \'etale.
Since $X \rightarrow C$ is a principal bundle, we see that after \'etale base change
$X \simeq C \times A$.

{\em 2nd case. The general  $\varphi$-fibre is algebraic.}

By Lemma \ref{lemmarelativedimensiontwo} there exists a normal projective surface $Z$ 
that admits an elliptic fibration \holom{\psi}{Z}{C} and satisfies the following properties:
\begin{itemize}
\item there exists an almost holomorphic fibration \merom{\tau}{X}{Z} such that $\varphi=\psi \circ \tau$.
\item there exists an endomorphism \holom{\overline{f}}{Z}{Z} of degree $d$ that commutes with $\psi$ and such that
$\overline{f} \circ \tau=\tau \circ f$. 
\end{itemize}

By \cite[Thm.6.2.1]{Nak08} there exists a finite 
base change $C' \rightarrow C$ such that the normalisation of $Z \times_C C'$ is isomorphic to $C' \times E$, where $E$ is an elliptic curve.
Note furthermore that by \cite[Lemma 6.2.5]{Nak08} the base change is \'etale over the locus where the fibres
of $\psi$ are reduced, in fact the ramification order in $c \in C$ equals the multiplicity of the fibre $Z_c$. 
Since $\varphi = \psi \circ \tau$ this shows that the ramification order in $c$ divides the multiplicity of the fibre $X_c$.
Thus by a classical argument the map from the normalisation of $X \times_C C'$ onto $X$ is \'etale.
As in the proof of Theorem \ref{theoremkappanminustwoA} we can use the universal property of the fibre product
and a commutative diagram to see that (up to \'etale cover) we can 
suppose that $Z \simeq C \times E$ where $E$ in elliptic curve.
By \cite[Lemma 2.25]{FN07} we may suppose (up to making an \'etale base change) that
the endomorphism $\overline{f}$ is of the form $id_C \times g_E$, where $\holom{g_E}{E}{E}$
is an endomorphism of degree $d$.

Since $E$ is 
an elliptic curve, the meromorphic fibration \merom{p_E \circ \tau}{X}{E} extends
to a holomorphic map such that we have a commutative diagram
$$
 \xymatrix{ 
X \ar @/_2pc/[dd]_{p_E \circ \tau} \ar @{-->}[d]^{\tau} \ar[r]^{f} & X  \ar @{-->}[d]^{\tau}  \ar @/^2pc/[dd]^{p_E \circ \tau} \\
C \times E  \ar[d]^{p_E} \ar[r]^{\overline{f}} & C \times E \ar[d]^{p_E} \\
E \ar[r]^{g_E} & E
}.
$$
Thus by Proposition \ref{propositionspecialfibres} the fibration $p_E \circ \tau$ is smooth and by adjunction
its fibres $X_e$ are surfaces with Kodaira dimension one that have a natural elliptic fibration \holom{\varphi|_{X_e}}{X_e}{C}. 
Since the general fibre of $\varphi$ is algebraic, an easy
application of Campana's theorem shows that a general fibre of $X_e$ has algebraic dimension one.
Thus the relative algebraic reduction of $X \rightarrow E$ is holomorphic and identifies to $\tau$. Thus $\tau$
is holomorphic and since $\overline{f}$ is \'etale of degree at least two, Proposition \ref{propositionspecialfibres}
implies that it is an equidimensional elliptic fibration whose singular locus is a disjoint union of $c_i \times E$.

Let now $X_c$ be a very general $\varphi$-fibre. Then $X_c$ is covered by an algebraic torus and is
an elliptic bundle $\tau_c: X_c \rightarrow c \times E \simeq E$. Thus by Poincar\'e's irreducibility theorem there exists an \'etale cover
of $X_c$ by $E'_c \times E$. In particular we get a family of elliptic curves in $X_c$ surjecting onto $E$. 
By countability of the components of the cycle space, we may suppose that the family of curves
deforms with $X_c$ for $c \in C$ very general. Let $Y \rightarrow C$ be the normalisation of the component of $\chow{X/C}$ parametrising 
these curves, and denote by $\tilde{X}$ the normalisation of the universal family. 
Denoting by \holom{p}{\tilde{X}}{X} and \holom{q}{\tilde{X}}{Y} the natural maps, we 
we have commutative diagram
$$
\xymatrix{ 
\tilde{X}  \ar[r]^{p} \ar[d]^{q} & X  \ar[d]^{\varphi} \ar[r]^{p_E \circ \tau} & E  
\\
Y \ar[r] & C
}.
$$
A general member of the family $\tilde{X}_y$ is an elliptic curve surjecting onto $E$. 
The \'etale base change  $\tilde{X}_y \rightarrow E$ induces an \'etale covering of $X \times_E \tilde{X}_y \rightarrow X$,
so up to replacing $X$ by an \'etale cover we can suppose without loss of generality that a general $\tilde{X}_y$
is a $p_E \circ \tau$-section and $p$ is birational.
We claim that $p$ is an isomorphism. Assuming this for the time being, let us show how to conclude:
for simplicity of notation, identify $X$ and $\tilde{X}$.
Let $X_e$ be a general $p_E \circ \tau$-fibre, then $X_e \rightarrow Y$ is surjective and
\'etale in codimension one. Thus the induced morphism $X \times_T X_e \rightarrow X$ is \'etale in codimension one,
so \'etale since $X$ is smooth. Thus up to replacing $X$ by an \'etale cover the $p_E \circ \tau$-fibres $X_e$
are $q$-sections, in particular all the fibres are isomorphic to $Y$. 
Thus $q \times (p_E \circ \tau): X \rightarrow Y \times E$ is an isomorphism.

{\em Proof of the claim.}
Since $X$ is smooth and $\tilde{X}$ is normal, it is sufficient that for every $x \in X$ there exist at most
finitely many $\tilde{X_y}$ passing through $x$. We will show this property ``fibrewise'':
let first $c \in C$ be a point that is not in the $\varphi$-singular locus. Then the fibre $X_c$ is an elliptic bundle over $c \times E$
and the $\tilde{X}_y$ are a one-dimensional family of disjoint sections parametrised by $Y_c$. Thus $X_c$ is a torus isomorphic
to $Y_c \times E$ and the $\tilde{X}_y$ form a family of line bundles such that the intersection product in
$X_c$ equals zero.

Let now $c_0 \in C$ be a point that is in the $\varphi$-singular locus. 
Then for every point $y \in Y_{c_0}$ there exists a small analytic neighbourhood $\D$ of $c_0$ and a 
multisection $S \subset Y$ over $\Delta$ passing through $y$ such that $S \cap Y_{c_0}$ is a singleton. 
Let $D:=p(\fibre{q}{S})$ be the analytic subset of $X$ covered by the curves parametrised by $S$. Then $D$ is a surface in the threefold $X$, so it is locally principal.
In particular for every $c \in \Delta$, the intersection $X_c \cap D$ is a Cartier divisor. Since for $c \in \D$ general the self-intersection is zero and the $X_c$
vary in a flat family, we see that the curve $\tilde{X}_y$ is a Cartier divisor in $X_c$ with self-intersection zero. 
Thus if we show that every curve $\tilde{X}_y$ is irreducible, this implies that there are at most finitely many $\tilde{X}_y$
through a given point and we are done.

The irreducibility of $\tilde{X}_y$ can be seen as follows: by Kodaira's list of singular fibres of an elliptic fibration, the reduction of the
fibre $X_c \simeq \fibre{\tau}{c \times E}$ is either an abelian surface (in this case the fibres over $c \times E$ are multiple elliptic curves)
or the fibres of $X_c \rightarrow E$ are cycles of (maybe singular) rational curves. We will deal with the second case, the first case can be reduced to
the preceding case by a local base change.
Since $f$ commutes with $\varphi$, we have an \'etale endomorphism \holom{f_c}{X_c}{X_c} of degree $d$ and up to replacing $f$ by some power 
$f_c$ maps every irreducible component of $X_c$ onto itself.  Again by Kodaira's list the normalisation $\holom{\nu}{T}{X_c}$ is a disjoint union of
$\PP^1$-bundles $T_i$ over the elliptic curve $E$ and $f_c$ lifts to an \'etale endomorphism $f_c^i: T_i \rightarrow T_i$ of degree $d$. 
Since $f^i_c$ is \'etale, one sees easily (cf. \cite[Ch.2]{Nak02}) that $T_i$ does not contain any curves with negative self-intersection.
The pull-back $\nu^* \tilde{X}_y$ is an effective divisor with self-intersection $0$ that surjects onto $E$. The 
claim follows by an easy intersection calculus on $T_i$.
\end{proof}

\subsection{Proof of Theorem \ref{maintheoremA}} \label{subsectionnonminimal}

By Proposition \ref{propositionmmpnotuniruled}
there exists a finite sequence of smooth compact K\"ahler threefolds
\[
X=X_0 \stackrel{\mu=\mu_1}{\longrightarrow} X_1 \stackrel{\mu_2}{\longrightarrow} \ldots \stackrel{\mu_n}{\longrightarrow} X_n
\]
such that $X_n$ is described by Theorem \ref{theoremminimalmodels} and $X_{i-1}$ is the blow-up of $X_i$ along
an elliptic curve $C_i$ with numerically trivial normal bundle such that $\fibre{f_i}{C_i}=C_i$.
We will now use the list in Theorem \ref{theoremminimalmodels} to see what blow-ups are actually possible.

{\em 1st case. $\kappa(X)=0$.} In this case $X_n$ is a torus or a product $Y \times E$ where $Y$ is a K3 surface 
and $E$ is an elliptic curve. 

{\em a) $X_n$ is a torus.} 
We claim that if $X \neq X_n$, then $X_n$ is isogenous to $Y \times E$ where $Y$ is a torus and $E$ an elliptic curve:
by hypothesis $X_n$ contains the elliptic curve $E:=C_n$, so we have a quotient map $\holom{\varphi}{X_n}{X_n/E}$. 
Since $\fibre{f_n}{C_n}=C_n$ 
there exists an \'etale endomorphism \holom{g}{X_n/E}{X_n/E} such that $g \circ \varphi=\varphi \circ f_n$.
Since $\fibre{f_n}{C_n}=C_n$ the fibre $\fibre{g}{\varphi(C_n)}$ is a singleton, so $g$ is an automorphism.
Thus $X_n$ is decomposable by Proposition \ref{propositionellipticovertorus}, i.e.\ a product of an elliptic curve
and a two-dimensional torus.
In particular the algebraic dimension of $X_n$ is at least one and if it equals one, the unique decomposition possible is $Y \times E$.
If the algebraic dimension of $X_n$ equals two, the quotient $X_n/E$ has algebraic dimension one and the claim follows from 
Proposition \ref{propositiontorusoverelliptic}.

Thus $X_n \simeq Y \times E$ and up to making an \'etale base change we can suppose that $f$ is of the form
$(g,h)$ where $\holom{g}{Y}{Y}$ is an automorphism and $\holom{h}{E}{E}$ an endomorphism of degree $d$ (cf. \cite[Lemma 2.25]{FN07} which does
not use the projectiveness assumption).
Let $0 \times C_n$ be the copy of $E \subset X_n$ that is the center of the blow-up. 
Then the commutative diagram
$$
\xymatrix{ 
X_n  \ar[r]^{f} \ar[d]^{p_E} & X_n  \ar[d]^{p_E}   
\\
E \ar[r]^h & E
}
$$
and $\deg h>1$ immediately implies that $C_n$ is not contained in some $Y \times e$ for $e \in E$.
Thus $C_n$ surject onto $E$ and up to replacing by an \'etale cover, it is of the form $y \times E$ for some $y_n \in Y$.
Thus the blow-up of $X_n$ along $C_n$ is isomorphic to $Bl_y Y \times E$.
The restriction of $f_{n-1}$ to $Bl_y Y$ is an automorphism, so we obtain the statement by
induction.

{\em b) $X_n \simeq Y \times E$ with $Y$ a K3.} 
Since $h^1(Y, \sO_Y)=0$, the projection onto $E$ is the Albanese map of $X_n$.
Since a K3 does not admit endomorphisms of degree at least two, 
the induced endomorphism $\holom{h}{E}{E}$ has degree $d$.
Since the automorphism group of a K3 surface is discrete we see that 
the restriction of $f_n$ to $Y \times e$ does not depend on $e$, so
$f_n=(g,h)$, where $\holom{g}{Y}{Y}$ is an automorphism.
We can now argue as in Case a).

{\em 2nd case. $\kappa(X)=1$.}
According to Theorem \ref{theoremminimalmodels} we distinguish two cases.

{\em a) $X_n \simeq C \times A$ with $C$ a curve of general type and $A$ a torus of algebraic dimension at most one.}
We claim that in this case it is not possible to have $X \neq X_n$ and argue by contradiction.
The projection onto $C$ is the Iitaka fibration,
so up to replacing $f$ by some power we may suppose that
 the induced endomorphism on $C$ is the identity. 
Moreover since the algebraic dimension of
$A$ is at most one, any morphism from $C$ to $A$ is constant. This implies that
$f_n = (Id_C, g)$ with $g$ an endomorphism of $A$ of degree $d$. The elliptic curve $C_n$ does not map surjectively onto
$C$, so it is contained in some torus $c \times A$. In particular $A$ has algebraic dimension one and $C_n$ is a fibre of the algebraic
reduction $A \rightarrow T$. Since $\fibre{g}{C_n}=C_n$, the endomorphism induced by $g$ on $T$ is an automorphism.
Hence $A$ is algebraic by Proposition \ref{propositiontorusoverelliptic}, a contradiction.

{\em b) $X_n \simeq E \times S$
where $E$ is an elliptic curve and $S$ a surface with $a(S)=\kappa(S)=1$.}
Let \holom{\psi}{S}{C} be the algebraic reduction of $S$. Then
\holom{p_E \times \psi}{X}{E \times C} is the algebraic reduction of $X$,
so there exists an induced endomorphism $\holom{g}{E \times C}{E \times C}$.
In general the endomorphism $g$ does not preserve the projection onto $E$, but by
\cite[Lemma 2.25]{FN07} there exists an \'etale base change $E' \rightarrow E$ such that this is the case. 
Since all the $X_i$ are fibre spaces over $E$ we can suppose without loss of generality that
there exists an endomorphism \holom{g_E}{E}{E} such that $g_E \circ p_E=p_E \circ g$.
Thus we also have $g_E \circ p_E=p_E \circ f$ and the restriction of $f$ to some $e \times S$ 
gives an endomorphism $f_S$of $S$. Since $S$ has algebraic dimension one, Theorem \ref{theoremkappanminusone}
implies that $f_S$ is an automorphism. Thus $g_E$ has degree $d$ and we can conclude as in the first case.
$\square$

\section{Uniruled manifolds} \label{sectionuniruled}

The aim of this section is to prove Theorem \ref{maintheoremB}. We will start with the easy case where $f$ is \'etale,
for the case of ramified endomorphisms we make a case distinction in terms of the algebraic dimension. With increasing
algebraic dimension the proofs get more and more involved, this confirms our philosophy that non-algebraicity is
an obstruction to the existence of endomorphisms.
For the convenience of the reader we will repeat at the start of every subsection 
the part of  Theorem \ref{maintheoremB} that we are about to prove.

\subsection{Proof of Theorem \ref{maintheoremB}: \'etale endomorphisms}

\begin{theorem} \label{theoremuniruledetale}
Let $X$ be a  non-algebraic compact K\"ahler threefold which is not uniruled, 
and let \holom{f}{X}{X} be an {\'etale} endomorphism of degree $d>1$.
Then  $X$ is (up to \'etale cover) a 
projectivised bundle $\PP(E)$ over a non-algebraic torus $A$ and $c_1^2(E)=4c_2(E)$.
\end{theorem}

\begin{proof}[Proof of Theorem \ref{theoremuniruledetale}]

Since the endomorphism $f$ is \'etale, we have $\chi(X, \sO_X)=0$.
Since $X$ is non-algebraic and uniruled, we have $h^3(X, \sO_X)=0$
and $h^2(X, \sO_X) \geq 1$. Therefore $h^1(X, \sO_X) \geq 2$, and we denote by $\alpha: X \rightarrow Y$
the Stein factorisation of the Albanese map
$X \rightarrow Alb(X)$ whose general fibre has dimension at least one.

{\em 1st case. $Y$ is a surface.}
The general fibre of $\alpha$ is a rational curve, so $Y$ is not algebraic since otherwise $X$ would be algebraic.
By Proposition \ref{propositioncommutation} the endomorphism $f$ induces an endomorphism $g$ on $Y$
which has degree $d$ since $f$ is \'etale. Applying Theorem \ref{theoremkappanminusone} we obtain $\kappa(Y)=0$,
so $Y \rightarrow Alb(X)$ is surjective and \'etale. By the universal property of the Albanese torus, we obtain 
$Y=Alb(X)$. 

We will now show that $\alpha: X \rightarrow Y=Alb(X)$ is a $\PP^1$-bundle. By Proposition \ref{propositionspecialfibres}
the $\alpha$-singular locus $\Delta$ is empty or has pure dimension one and satisfies $\fibre{g}{\Delta}=\Delta$.
Since $Y$ is a non-algebraic two-dimensional torus, the components of $\Delta$ are elliptic curves
that are contracted by the algebraic reduction $Alb(X) \rightarrow E$. The endomorphism
$g$ induces an endomorphism $h$ on $E$ by Proposition \ref{propositionalgebraicreduction}.  
The condition $\fibre{g}{\Delta}=\Delta$ implies that $h$ is an isomorphism which contradicts $a(Alb(X)) \leq 1$
by Proposition \ref{propositiontorusoverelliptic}.

Up to making an \'etale base change, we can suppose that $X$ is a projectivised bundle $\PP(E)$.
The canonical bundle of a is trivial, so $K_X \simeq \sO_{\PP(E)}(-2) \otimes \varphi^* \det E$.
Since $f$ is \'etale, we have $K_X^3=0$. An elementary computation shows that $c_1^2(E)=4c_2(E)$.

{\em 2nd case. $Y$ is a curve.} We will prove that in this case $X$ is necessarily algebraic.
Since $h^1(Y, \sO_Y)=h^1(X, \sO_X) \geq 2$, the endomorphism induced on $Y$ is an automorphism
of finite order, so up to replacing $f$ by some power we may suppose that $g=Id_Y$.
Let $S \rightarrow C$ be the unique irreducible component of the relative cycle space $\chow{X/C}$
that parametrizes the rational curves in $X$.
The push-forward of cycles $f_*$ acts on $\chow{X/C}$ and since $f$ induces a meromorphic endomorphism
on the base of the rationally connected quotient $X \dashrightarrow S$ (Prop. \ref{propositioncommutation}).
Thus we obtain an endomorphism \holom{f_*}{S}{S} that commutes with $S \rightarrow C$.
Since $f$ is \'etale and the general fibre of $X \dashrightarrow S$ is a $\PP^1$, we see that
$f_*$ has degree $d$. The general fibre of $S \rightarrow C$ is an elliptic curve since otherwise $S$ is algebraic.
By Lemma \ref{lemmamultisection} we see that $S \rightarrow C$ has a multisection,
so $S$ is algebraic. Therefore $X$ is algebraic.
\end{proof}

\subsection{Proof of Theorem \ref{maintheoremB}: algebraic dimension zero}

\begin{theorem} \label{theoremzero}
Let $X$ be a compact K\"ahler threefold of algebraic dimension zero which is uniruled, 
and let \holom{f}{X}{X} be a ramified endomorphism of degree $d>1$.
Then $X$ is (up to \'etale cover) a projectivised bundle $\PP(E)$ over a torus $A$ of algebraic dimension zero,
$f$ induces an endomorphism on $A$ of degree at least two, and $E$ is a direct sum of line bundles.
\end{theorem}

If the algebraic dimension is zero, we have the following elementary, but useful lemma at our disposition.

\begin{lemma} \label{lemmafibrationzero}
Let $X$ be a compact K\"ahler threefold of algebraic dimension zero and let \holom{f}{X}{X} be an endomorphism of degree $d>1$. 
Suppose that $X$ admits a fibration \holom{\varphi}{X}{S} onto a normal surface $S$.
Then there exists an endomorphism $\holom{g}{S}{S}$ such that $g \circ \varphi=\varphi \circ g$.
\end{lemma}

\begin{proof}
By a theorem of Fischer and Forster \cite[Thm.]{FF79} there exist only finitely many divisors on $X$.
Since the images and preimages of divisors are divisors, we see that $f$ induces a bijective map on the set of divisors on $X$.
Since this set is finite, some iteration $f^k$ induces the identity. 

By the rigidity lemma it is sufficient to prove that any $\varphi$-fibre is mapped by $f$ into a $\varphi$-fibre.
By what precedes we may suppose that this is true for the divisorial fibre components.
Since $S$ contains only finitely many curves, it also holds for the general $\varphi$-fibres.
\end{proof}

\begin{proof}[Proof of Theorem \ref{theoremzero}]
The rationally connected quotient is an almost holomorphic map
\merom{\varphi}{X}{S} onto a compact K\"ahler surface of algebraic dimension zero.
By Lemma \ref{lemmaK3torus} we can replace $S$ by a singular K3 surface or a torus
without curves, so we may suppose that $\varphi$ is holomorphic. 
By Lemma \ref{lemmafibrationzero} there exists a holomorphic map \holom{g}{S}{S} such that
$\varphi \circ f=g \circ \varphi$.
Since $S$ does not contain any curves, the $\varphi$-singular locus is a finite union of points.

{\it 1st case. $g$ is an automorphism.}
Since the $\varphi$-singular locus is a finite union of points, the conditions of Proposition \ref{propositionamerik}
are satisfied (if $S$ is a singular K3 surface, we also use Corollary \ref{corollaryK3}). 
Thus  after finite base change $X$ is bimeromorphic  to $S \times \PP^1$. In particular it has algebraic dimension at least one, a contradiction.

{\it 2nd case. $g$ is not an automorphism.} 
We claim that $S$ is a torus: otherwise $S$ would be  a singular K3 surface and
by Proposition \ref{propositionkummer} there exists a Galois cover by a torus \holom{\nu}{A}{S} that is \'etale in codimension one
such that $g$ lifts to an endomorphism $\holom{g_A}{A}{A}$ of degree $\deg g$.
The following commutative diagram and the universal property of the fibre product show that $f$ lifts to an endomorphism
$\holom{f'}{X \times_S A}{X \times_S A}$.
$$
 \xymatrix{ 
X \times_S A \ar[rd] \ar[ddd]^{\varphi'} \ar @{.>}[rrr]_{f'}^{\exists \ \mbox{\small univ. prop. fibre prod.}} & &  &  X \times_S A  \ar[ld] \ar[ddd]^{\varphi'}
\\
& X \ar[d]^{\varphi} \ar[r]^{f} & X  \ar[d]^{\varphi} &  
\\
& S  \ar[r]^{g} & S & 
\\
A \ar[ru]^\nu \ar[rrr]^{g_A} & & & A \ar[lu]_\nu
}
$$
Since the endomorphism $g_A$ is \'etale and $A$ does not contain any curves, 
the morphism $\varphi'$ is smooth by Proposition \ref{propositionspecialfibres}.
Thus the reduction of every $\varphi$-fibre is isomorphic to $\PP^1$. Suppose that there exists a non-reduced fibre $\fibre{\varphi}{s}$: 
then $-K_X \cdot \fibre{\varphi}{s}=2$ implies that the fibre is a double $\PP^1$ and the reduced fibre $l$ satisfies $-K_X \cdot l=1$. 
By a theorem of Ein-Koll\'ar \cite[Thm.5.3]{Kol91} the rational curve $l$ deforms in a one-dimensional family that covers a surface $D$ in $X$. Since $\varphi$
is equidimensional,  $\varphi(D)$ is a curve in $S$, a contradiction.
Thus we see that $\varphi$ is a smooth map and hence $S$ is a smooth K\"ahler surface of algebraic dimension zero
that admits an endomorphism of degree at least two.
This shows that $S$ is a torus.

Up to replacing $X$ by an \'etale cover, we can suppose that $X=\PP(E)$ where $E$ is a rank two vector bundle
which splits by Lemma \ref{lemmaprojectivebundleovertoruszero} below.
\end{proof}

The same proof shows the following statement.

\begin{proposition} \label{propositionbasezero}
Let $X$ be a compact K\"ahler threefold, and let \holom{f}{X}{X} be an ramified endomorphism of degree $d>1$.
Suppose that the rationally connected quotient is a fibration $\holom{\varphi}{X}{S}$ onto a normal K\"ahler surface without curves.
Suppose furthermore that the induced endomorphism \holom{g}{S}{S} has degree at least two.
Then $S$ is a torus and $X$ a $\PP^1$-bundle over $S$.
\end{proposition}

\begin{lemma} \label{lemmaprojectivebundleovertoruszero}
Let $A$ be a two-dimensional torus of algebraic dimension zero, and let $\holom{\varphi}{X=\PP(E)}{A}$ be 
the projectivisation of a rank two vector bundle $E$. Suppose that there exists a ramified endomorphism
\holom{f}{X}{X} and an \'etale endomorphism \holom{g}{A}{A} of degree at least two such that $g \circ \varphi=\varphi \circ f$.
If the algebraic dimension of $X$ is zero, then $E$ is (up to finite \'etale base change) a direct sum of line bundles.
\end{lemma}

\begin{remark*}
The statement of the lemma as well as the  techniques in the proof are similar to \cite[Thm.2]{Ame03}.
Nevertheless there is an important difference due to the condition on $A$: our conclusion holds after {\em \'etale} base change. 
\end{remark*}

\begin{proof}
Since $X$ has algebraic dimension zero, there exist only finitely many divisors on $X$. Thus up to replacing $f$ by some multiple we can suppose
that $\fibre{f}{D}=D$ for every irreducible effective divisor on $X$. Since $f$ is ramified, there exists an effective divisor $D \subset R$
and we denote by \holom{f_D}{D}{D} the restriction of $f$ to $D$. The torus $A$ does not contain any curve, so the holomorphic map
\holom{\varphi_D}{D}{A} is surjective and not ramified. Since $g \circ \varphi_D=\varphi_D \circ f_D$ and $g$ \'etale of degree at least two,
Proposition \ref{propositionspecialfibres} shows that $\varphi_D$ does not have any higher-dimensional fibres.
Thus $\varphi_D$ is \'etale and $D$ is a torus of algebraic dimension zero. Making an \'etale base change $D \rightarrow A$, we can suppose that
$D$ is a $\varphi$-section. 
Moreover the following commutative diagram and the universal property of the fibre product show that $f$ lifts to an endomorphism
$\holom{f'}{X \times_A D}{X \times_A D}$.
$$
 \xymatrix{ 
X \times_A D \ar[rd] \ar[ddd] \ar @{.>}[rrr]_{f'}^{\exists \ \mbox{\small univ. prop. fibre prod.}} & &  &  X \times_A D  \ar[ld] \ar[ddd]
\\
& X \ar[d]^{\varphi} \ar[r]^{f} & X  \ar[d]^{\varphi} &  
\\
& A  \ar[r]^{g} & A & 
\\
D \ar[ru]^{\varphi_D} \ar[rrr]^{f_D} & & & D \ar[lu]_{\varphi_D}
}
$$
Thus up to finite \'etale base change the ramification divisor $R$ has one irreducible component $D_1$ that is a $\varphi$-section.
Since any endomorphism of $\PP^1$ of degree at least two ramifies in at least two points, this implies that $R$ has at least another irreducible component $D_2$.
Arguing as before we see that $D_2$ is an \'etale cover of $A$. Moreover the intersection with $D_1$ is empty since $D_1$ does not contain any curves.
Up to making a second \'etale base change $D_2 \rightarrow A$ we can suppose that $\varphi$ has two disjoint sections.
Thus $E$ is a direct sum of line bundles.
\end{proof}

\subsection{Proof of Theorem \ref{maintheoremB}: algebraic dimension one}

\begin{theorem} \label{theoremone}
Let $X$ be a compact K\"ahler threefold of algebraic dimension one which is uniruled, 
and let \holom{f}{X}{X} be a ramified endomorphism of degree $d>1$.
Then (up to \'etale cover) one of the following holds:
\begin{enumerate}
\item $X$ is a product $S \times \PP^1$, where $S$ is  a compact K\"ahler surface of algebraic dimension zero
and $f$ induces an automorphism on $S$.
\item $X$ is a projectivised bundle $\PP(E)$ over a torus $A$ of algebraic dimension at most one
such that the induced morphism on $A$ has degree at least two and $E$ is a direct sum of line bundles.
\end{enumerate}
\end{theorem}

\begin{proof}[Proof of Theorem \ref{theoremone}]
The rationally connected quotient is 
an almost holomorphic map \merom{\varphi}{X}{S} onto a compact K\"ahler surface of algebraic dimension at most one.

\smallskip

{\em 1st case. $S$ has algebraic dimension one}

\smallskip

Up to replacing $S$ by a bimeromorphic model, we may suppose that $S$ is a relatively minimal elliptic surface \holom{\psi}{S}{C}.
Since $a(S)=1$ there does
not exist any curve in $S$ that surjects onto $C$. 
The general fibre of  $\merom{\psi \circ \varphi}{X}{C}$ is uniruled, hence algebraic, so $\psi \circ \varphi$ extends to a holomorphic map
by Corollary  \ref{corollarygeneralfibrealgebraic}.
Furthermore we know by \cite[Cor.7.3]{CP00} that the general $\psi \circ \varphi$-fibre is isomorphic to $\PP(\sO_E \oplus L) \rightarrow E$
where $E$ is an elliptic curve and $L$ is a numerically trivial line bundle that is not torsion.

By Proposition \ref{propositioncommutation} there exists a meromorphic endomorphism \merom{g}{S}{S} such that 
$g \circ \varphi = \varphi \circ f$. Thus by Proposition \ref{propositionalgebraicreduction} there exists an endomorphism \holom{g_C}{C}{C}
such that $g_C \circ \psi = \psi \circ g$. Thus we get a commutative diagram
$$
 \xymatrix{ 
X \ar @/_2pc/[dd]_{\psi \circ \varphi} \ar @{-->}[d]^{\varphi} \ar[r]^{f} & X  \ar @{-->}[d]^{\varphi}  \ar @/^2pc/[dd]^{\psi \circ \varphi} \\
S  \ar[d]^{\psi} \ar @{-->}[r]^{g} & S \ar[d]^{\psi} \\
C  \ar[r]^{g_C} & C
}
$$
and make another case distinction.

{\it a) $g_C$ is an automorphism.} We will show that this case does not exist:
the restriction of $f$ to the general fibre $X_c \simeq \PP(\sO_E \oplus L)$ 
gives an endomorphism \holom{f_c}{X_c}{X_c} degree $d>1$.  Denote by \holom{\pi}{X_c}{E} 
the canonical projection, then there exists an endomorphism \holom{g_E}{E}{E}
such that $g_E \circ \pi = \pi \circ f_c$. Note that $g_E$ can't be an automorphism: otherwise
Amerik's theorem \cite[Thm.1]{Ame03} would imply that $X_c$ is a product after \'etale cover, so $L$ is a torsion
line bundle.
We will now use the Mori program: since $X$ is uniruled, it admits a contraction.

Suppose first that $X$ admits a fibre type contraction. 
Since $X$ contains only one covering family of rational curves and
the rationally connected quotient is only defined up to birational equivalence, we can suppose
that the contraction is the rationally connected quotient \holom{\varphi}{X}{S}.
Since $\varphi$ is flat, the rigidity lemma implies that $g$ extends to a {\em holomorphic}
endomorphism \holom{g}{S}{S}. Since the restriction of $g$ to a general $\psi$-fibre is $g_E$,
we see that $g$ has degree at least two. Thus $S$ is a torus \cite{FN05} and admits an endomorphism
$g$ of degree at least two such that  $g_C \circ \psi=\psi \circ g$. Since $g_C$ is an automorphism,
Proposition \ref{propositiontorusoverelliptic}  implies that $S$ is isogenous to a product of elliptic curves,
a contradiction to $a(S)=1$.

If $X$ admits a contraction of birational type \holom{\mu}{X}{X'}, we argue using
Proposition \ref{propositionkaehlerpackage}.
First observe that the composed map $\psi \circ \varphi: X \to C$ is locally projective (blow up to make $\varphi$ a (locally projective) morphism, then a 
priori
$\psi \circ \varphi$ might only locally Moishezon, but \cite[Thm.10.1]{CP04} gives local projectivity). 
Then change possibly $\mu$ to get local projectivity of the induced map $\tau': X \to C$.
Then we get inductively a sequence of a compact manifolds bimeromorphically equivalent to a K\"ahler manifold             
\[
X=X_0 \rightarrow X_1 \rightarrow \ldots \rightarrow X_n
\]
such that $X_n$ admits a fibre type contraction 
and $g_C$ is an automorphism. But we have just seen that such a $X_n$ does not exist (the K\"ahler property is not needed there).

{\it b) $g_C$ is not an automorphism.} 
In this case the curve $C$ is elliptic or $\PP^1$. 
We can't have 
$C \simeq \PP^1$, since otherwise $S$ is algebraic by Proposition \ref{propositionfavre}.
Thus $C$ is elliptic and the endomorphism $g_C$ is \'etale.
It follows by Lemma \ref{lemmafibrationcurve} that $\psi \circ \varphi$ is a submersion, 
thus all the fibres are uniruled and have $b_1(X_c)=b_2(X_c)=2$. This easily implies that all the fibres are 
$\PP^1$-bundles over elliptic curves. Thus $\varphi$ extends to a holomorphic map and $X$ is a $\PP^1$-bundle over $S$
which is an elliptic bundle over $S$. Since $S$ is K\"ahler but not algebraic, it follows from  \cite[V.5.B)]{BHPV04}
that $S$ is a torus of algebraic dimension one. Thus the meromorphic endomorphism  \merom{g}{S}{S}
also extends to a holomorphic map of degree at least $\deg g_C$.
Up to making an \'etale base change, we can suppose
$X \simeq \PP(E)$, where $E$ is a rank two vector bundle which splits by Lemma \ref{lemmaprojectivebundleovertorusone} below.

\smallskip

{\em 2nd case. $S$ has algebraic dimension zero}

\smallskip

In this case $X$ is bimeromorphic to $\PP^1 \times S$ where $S$ is a compact K\"ahler surface of algebraic dimension zero
\cite[Cor.7.6]{CP00}. The general fibre of the algebraic reduction is bimeromorphic to $S$, so by Proposition \ref{propositionexistencefibration}
we have a holomorphic algebraic reduction $\holom{r}{X}{\PP^1}$.
On the other hand, we know by Lemma \ref{lemmaK3torus} that, up to replacing $S$ by some bimeromorphic normal model, we may 
suppose that $S$ does not contain any curves. Thus we may suppose that the rationally connected quotient is a holomorphic map \holom{\varphi}{X}{S},
so we get a {\em holomorphic} map \holom{r \times \varphi}{X}{\PP^1 \times S} of degree one.
In particular we have 
$$r^* \sO_{\PP^1}(1) \cdot f=1,$$ 
where $f$ is a general $\varphi$-fibre.
We know by Proposition \ref{propositioncommutation} that there exists an endomorphism \holom{g_S}{S}{S} such that 
$g_S \circ \varphi = \varphi \circ f$. Moreover  by Proposition \ref{propositionalgebraicreduction} there exists an endomorphism \holom{g_{\PP^1}}{\PP^1}{\PP^1}
such that $g_{\PP^1} \circ r = r \circ f$. Together they induce an endomorphism
\holom{g_{\PP^1} \times g_S}{\PP^1 \times S}{\PP^1 \times S} such that we have a commutative diagram
$$
 \xymatrix{ 
X \ar[d]^{r \times \varphi} \ar[r]^{f} & X  \ar[d]^{r \times \varphi}  \\
\PP^1 \times S  \ar[r]^{g_{\PP^1} \times g_S} & \PP^1 \times S}
$$

We distinguish two cases:

{\it a) $g_S$ is not an automorphism.}
By Proposition \ref{propositionbasezero}  we know
that $X \rightarrow S$ is a $\PP^1$-bundle and $S$ is a torus. 
Since $r^* \sO_{\PP^1}(1) \cdot f=1$ the fibres of the algebraic reduction are $\varphi$-sections.
Thus $r \times \varphi$ is an isomorphism onto $\PP^1 \times S$ and we are in the second Case of Theorem \ref{theoremone} (in this
case $E \simeq \sO_S^{\oplus 2}$).

\medskip

{\it b) $g_S$ is an automorphism.}
Since $r \times \varphi$  has degree one, this implies that $g_{\PP^1}$ has degree $d$.
Let $T_2 \subset \PP^1 \times S$ be the set such that the $r \times \varphi$-fibre has dimension two.
Then $T_2$ is finite and by Proposition \ref{propositionspecialfibres}, the endomorphism $g_{\PP^1} \times g_S$
is totally ramified in every point of $T_2$.
Since $g_S$ is an automorphism, this shows that $g_{\PP^1}$ has ramification order $d$ in every point of $p_{\PP^1}(T_2)$.

We will now use Mori theory to discuss the structure of $X$: since $X$ is uniruled it admits at least one contraction.

{\it b1) $X$ admits a fibre type contraction \holom{\psi}{X}{S'}.} 
Then $X$ is a $\PP^1$- or conic bundle over the smooth surface $S'$ (Theorem \ref{theoremcontraction}). 
Since $X$ contains only one covering family of rational curves, the general $\psi$-fibre is a general $\varphi$-fibre. The fibration $\psi$ is flat, so the rigidity lemma
shows that we have a factorisation \holom{\tau}{S'}{S} which is actually a birational map. 
Since $\rho(X/S')=1$ and $r^* \sO_{\PP^1}(1) \cdot f=1$,
we see that $\psi$ is not a conic bundle. 
The endomorphism $g_S$ is an automorphism, so the endomorphism $g_{S'}$ induced by $f$ on $S'$ is an automorphism.
We satisfy the conditions of Proposition \ref{propositionamerik}, so there exists an \'etale cover of $S'$ such that $X$ becomes
a product after base change. Thus we are in the first case of Theorem \ref{theoremone}.
Note furthermore that $f=(g_S,h)$ since the space of endomorphism of $\PP^1$ is affine \cite[Lemma 1.2]{Ame03}.

{\it b2) $X$ admits a birational contraction \holom{\psi}{X}{X'}.}
Denote by $E$ the exceptional locus. Since $S$ does not contain
any curve, $E$ is contained in a higher-dimensional fibre of $\varphi$. Since there are only finitely many higher-dimensional fibres,
we can suppose (up to replacing $f$
by $f^k$) that $\fibre{f}{E}=E$. In particular we get an endomorphism \holom{f'}{X'}{X'} of degree $d$ such that $f' \circ \psi=\psi \circ f$.
We claim that $\psi$ is of type 3.) in Theorem \ref{theoremcontraction}, i.e.\ it contracts a divisor onto a smooth curve.

{\it Proof of the claim:} We argue by contradiction and suppose that $\psi$ contracts a divisor $E$ onto a point. Since all the curves contained in $E$ 
are numerically equivalent and $\PP^1 \times S$ does not contain any ruled surface (recall that $S$ contains no curve), we see that
$r \times \varphi$ maps $E$ onto a point $c$ contained in $T_2$. By what precedes, the endomorphism
$g_{\PP^1}$ has ramification order $d$ in $p_{\PP^1}(c)$. Denote by $X_c=r^* c$ the $r$-fibre over $c$.
Since $g_{\PP^1} \circ r = r \circ f$ and $g_{\PP^1}^* c = d c$, we have
$$
d X_c = f^* X_c.
$$
Thus the ramification order of $f$ along $E \subset X_c$ equals $d$ and $E^3=N_{E/X}^2\neq 0$ by Theorem \ref{theoremcontraction}.
Since $\fibre{f}{E}=E$, this contradicts Corollary \ref{corollarycontractionuniruled}. $\square$

Thus $X'$ is smooth and
since $S$ does not contain any curves, there exists a holomorphic map $\varphi': X' \rightarrow S$
such that $\varphi=\varphi' \circ \psi$.
Note that we also obtain a morphism
$$ t: X' \to \PN_1 \times S.$$
Instead of applying Proposition \ref{propositionkaehlerpackage}, which is tedious, we argue as follows. We apply directly \cite{Nak87} to the {\it projective}
morphism 
$$ r \times \varphi: X \to \PN_1 \times S $$
and obtain a possibly new birational morphism $\nu: X \to X'$ and a factorisation $X' \to \PN_1 \times S$ which is again projective. Since $\PN_1 \times S$ is
K\"ahler, so is $X'$.    
Arguing by induction, we get a sequence of birational contractions
$$
X=X_0 \rightarrow X_1 \rightarrow \ldots \rightarrow X_n
$$
such that $X_n$ is compact K\"ahler, admits an endomorphism  \holom{f_n}{X_n}{X_n} of degree $d$ and a fibre type contraction.  
By Case a) there exists an \'etale cover of  $X_n$  by $S_n \times \PP^1$ 
and 
$f_n=(g_n,h_n)$ where $g_n$ is an automorphism of $S_n$ and $h_n: \PP^1 \rightarrow \PP^1$ an endomorphism
of degree $d$.
Since 
$$
\pi_1(X) \simeq \pi_1(X_1) \simeq \ldots \simeq \pi_1(X_n)
$$ 
we can suppose that  $X_n=S_n \times \PP^1$.
The manifold $X_{n-1}$ is obtained by blowing up a curve $C_n$ such that $\fibre{f_n}{C_n}=C_n$. It is elementary 
to see that such a curve is necessarily of the form $s_n \times \PP^1$ where $s_n \in S_n$ is a point. Thus
$X_{n-1} \simeq Bl_{s_n} S_n \times \PP^1$ and we conclude inductively that we are in the first case of Theorem \ref{theoremone}.
\end{proof}

\begin{lemma} \label{lemmaprojectivebundleovertorusone}
Let $A$ be a two-dimensional torus of algebraic dimension one, and let $\holom{\varphi}{X=\PP(E)}{A}$ be 
the projectivisation of a rank two vector bundle $E$. Suppose that there exists a ramified endomorphism
\holom{f}{X}{X} and an \'etale endomorphism \holom{g}{A}{A} of degree at least two such that $g \circ \varphi=\varphi \circ f$.
If the algebraic dimension of $X$ is one, then $E$ is (up to \'etale cover) a direct sum of line bundles.
\end{lemma}

\begin{proof}
Let \holom{\psi}{A}{E} be the algebraic reduction of the torus $A$, and denote by \holom{g_E}{E}{E} the endomorphism induced by $g$ on $E$.
Since $X$ has algebraic dimension one, $\psi \circ \varphi$ is the algebraic reduction of $X$. 
Since $f$ is ramified and $g$ \'etale any effective divisor $D \subset R$ maps surjectively onto $A$,
in particular it is not a polar divisor. Since there are only finitely many non-polar divisors \cite[Thm.]{FF79},
we can suppose (up to replacing $f$ by some power) that $\fibre{f}{D}=D$ for every irreducible component of $D \subset R$ effective divisor on $X$.

Denote by \holom{f_D}{D}{D} the restriction of $f$ to some  $D \subset R$. 
Since $g \circ \varphi_D=\varphi_D \circ f_D$ and $g$ is \'etale of degree at least two,
Proposition \ref{propositionspecialfibres} shows that $\varphi_D$ does not have any higher-dimensional fibres.
If $\varphi_D$ is ramified, the branch locus $\Delta \subset A$ is a finite union of curves such that $\fibre{g}{\Delta}=\Delta$.
Yet this implies that $g_E$ is an automorphism, so $A$ would be algebraic by Proposition \ref{propositiontorusoverelliptic}.
Thus $\varphi_D$ is \'etale and $D$ a torus. Arguing as in the proof of Lemma \ref{lemmaprojectivebundleovertoruszero}
we see that up to \'etale base change, the ramification $R$ has at least two irreducible components: one component 
$D_1$ that is a $\varphi$-section and a second $D_2$ that is an \'etale cover of $A$.  
Moreover the intersection with is empty since
$D_1 \cap D_2$ would be a bunch of curves such that $\fibre{f_{D_1}}{D_1 \cap D_2} = D_1 \cap D_2$. 
Once more Proposition \ref{propositiontorusoverelliptic} would then imply that $D_1 \simeq A$ is algebraic.
Thus up to making a second \'etale base change $D_2 \rightarrow A$, we can suppose that $\varphi$ has two disjoint sections.
Thus $E$ is a direct sum of line bundles.
\end{proof}

\subsection{Proof of Theorem \ref{maintheoremB}: algebraic dimension two}

\begin{theorem} \label{theoremtwo}
Let $X$ be a compact K\"ahler threefold of algebraic dimension two which is uniruled, 
and let \holom{f}{X}{X} be a ramified endomorphism of degree $d>1$.
Then (up to \'etale cover) one of the following holds:
\begin{enumerate}
\item $X$ is a prodcut $Y \times \PP^1$ where $Y$ is a surface of algebraic dimension one 
and $f$ induces an automorphism on $Y$ or
\item $X$ is a projectivised bundle $\PP(E)$ over a torus $A$ of algebraic dimension one and $f$ induces
an endomorphism $g$ of degree at least two on $A$.
\end{enumerate}
\end{theorem}

\begin{remarks}
1. It is obvious that the list is effective: in the second case take $X:=A \times \PP^1$ with $A$ a torus of algebraic dimension one.
The following example shows that in the second case $X$ is in general not a product after \'etale cover.

Let $E$ be an elliptic curve without complex multiplication, and let $L$ be a numerically trivial line bundle on $E$ that is not torsion.
By \cite[Prop.5]{Nak02} there exists a ramified endomorphism
\holom{h}{\PP(\sO_E \oplus L)}{\PP(\sO_E \oplus L)} 
preservering the fibration such that $\PP(\sO_E \oplus L) \rightarrow E$
and inducing an endomorphism
\holom{g_E}{E}{E} of degree at least two.
Since $End(E)\simeq \Z$ it is clear that $g_E$ is of the form $x \mapsto n x + e$ for some $e \in E$.
Let now $\psi: A \rightarrow E$ be a torus of algebraic dimension one,
and let  \holom{g}{A}{A} be the endomorphism $x \mapsto n x + a$ where $a \in \fibre{\psi}{e}$.
Then we obviously have $g_E \circ \psi=\psi \circ g$.

The fibre product $X:=\PP(\sO_E \oplus L) \times_E A$ is a  $\PP^1$-bundle
over $A$ that is not a product, even after \'etale cover. An easy diagram chase shows that $X$ admits
a ramified endomorphism $f$ preserving the $\PP^1$-bundle structure and inducing the endomorphism $g$ on $A$.

2. One might also ask for a more precise classification of the vector bundle $E$ in the second case of the theorem: if
$A \rightarrow B$ is the algebraic reduction of $A$, it is not hard to see $E|_{A_b} \simeq L_b \oplus L_b \otimes T_b$,
where $T_b$ is a torsion line bundle on $A_b$. Since $L_b$ and $T_b$ might deform with $b \in B$, a further discussion
is quite tedious and we leave it as an exercise to the interested reader.
\end{remarks}

\begin{proof}[Proof of Theorem \ref{theoremtwo}]           

The rationally connected quotient is 
an almost holomorphic map \merom{\varphi}{X}{S} onto a compact K\"ahler surface of algebraic dimension one \cite[Thm.9.1]{CP00}.

Up to replacing $S$ by a bimeromorphic model, we may suppose that $S$ is a relatively minimal elliptic surface \holom{\psi}{S}{C}.
Since $a(S)=1$ there does not exist any curve in $S$ that surjects onto $C$. 
The general fibre of  $\merom{\psi \circ \varphi}{X}{C}$ is uniruled, hence algebraic, so $\psi \circ \varphi$ extends to a holomorphic map
by Corollary  \ref{corollarygeneralfibrealgebraic}.

By Proposition \ref{propositioncommutation} there exists a meromorphic endomorphism \merom{g}{S}{S} such that 
$g \circ \varphi = \varphi \circ f$. Thus by Proposition \ref{propositionalgebraicreduction} there exists an endomorphism \holom{g_C}{C}{C}
such that $g_C \circ \psi = \psi \circ g$. Thus we get a commutative diagram
$$
 \xymatrix{ 
X \ar @/_2pc/[dd]_{\psi \circ \varphi} \ar @{-->}[d]^{\varphi} \ar[r]^{f} & X  \ar @{-->}[d]^{\varphi}  \ar @/^2pc/[dd]^{\psi \circ \varphi} \\
S  \ar[d]^{\psi} \ar @{-->}[r]^{g} & S \ar[d]^{\psi} \\
C  \ar[r]^{g_C} & C
}
$$

{\it 1st case. $g_C$ is not an automorphism.} 
In this case the curve $C$ is elliptic or $\PP^1$. We can't have $C \simeq \PP^1$, since otherwise $S$ is algebraic by Proposition \ref{propositionfavre}.
Thus $C$ is elliptic and the endomorphism $g_C$ is \'etale.
It follows by Lemma \ref{lemmafibrationcurve} that $\psi \circ \varphi$ is a submersion. 
The general fibre $X_c$ is uniruled and the base of its rationally connected quotient is an elliptic curve $S_c=\fibre{\psi}{c}$,
so $b_1(X_c)=2$. Since $\psi \circ \varphi$ is smooth, this implies that
the rationally connected connected quotient of all the fibres  is a holomorphic fibration 
over an elliptic curve. Thus $\varphi$ extends to a holomorphic map and $S$
is an elliptic bundle over the elliptic curve $C$. Since $S$ is K\"ahler but not algebraic, it follows from  \cite[V.5.B)]{BHPV04}
that $S$ is a torus of algebraic dimension one. Thus the meromorphic endomorphism  \merom{g}{S}{S}
extends to a holomorphic \'etale map of degree at least $\deg g_C>1$. 
Thus by Proposition \ref{propositionspecialfibres} the $\varphi$-singular locus $\Delta$ is empty or a finite union of curves
such that $\fibre{g}{\Delta}=\Delta$. Since $S$ has algebraic dimension one, the irreducible components of $\Delta$ are
$\psi$-fibres. Thus $\fibre{g}{\Delta}=\Delta$ implies $\fibre{g_C}{\psi(\Delta)}=\psi(\Delta)$. Since $g_C$ is \'etale, we see that $\psi(\Delta)$ is empty.
Thus $X$ is a $\PP^1$-bundle over the torus $S$ and (up to another \'etale base change)
we are in the second case of Theorem \ref{theoremtwo}.

{\it 2nd case. $g_C$ is an automorphism.}
We will use Mori theory to discuss the structure of $X$: since $X$ is uniruled it admits at least one contraction.

{\it a) $X$ admits a fibre type contraction.}
Since $X$ contains only one covering family of rational curves and
the rationally connected quotient is only defined up to bimeromorphic equivalence, we can suppose
that the contraction is the rationally connected quotient \holom{\varphi}{X}{S}.
Then $X$ is a $\PP^1$- or conic bundle over the smooth surface $S$ (Theorem \ref{theoremcontraction})
which is an elliptic surface over $C$.
Since $\varphi$ is flat, the rigidity lemma implies that $g$ extends to a {\em holomorphic}
endomorphism \holom{g}{S}{S}.

{\it a1) $g$ is not an automorphism.}
In this case $S$ is a torus and admits an endomorphism
$g$ of degree at least two such that  $g_C \circ \psi=\psi \circ g$. Since $g_C$ is an automorphism,
Proposition \ref{propositiontorusoverelliptic}  implies that $S$ is isogenous to a product of elliptic curves,
a contradiction to $a(S)=1$.

{\it a2) $g$ is an automorphism.}
.We claim that in this case the variety $X$ admits a {\em holomorphic} elliptic fibration
\holom{r}{X}{Y} such that $a(Y)=2$, i.e.\ the algebraic reduction of $X$ can be taken holomorphic.

{\em Proof of the claim.} Let $M \subset S$ be an irreducible curve. By Kodaira's classification
of fibres of elliptic fibrations, the normalisation $\tilde{M} \rightarrow M$ is an elliptic or rational curve.
Set 
$$X_M :=\fibre{\varphi}{M},$$ 
then the general fibre of $\holom{\varphi_M}{X_M}{M}$ is a $\PP^1$
or two $\PP^1$'s meeting transversally in a point. Moreover (up to replacing $f$ by some power) 
the restriction of $f$ to $X_M$ gives an endomorphism
$$\holom{f_M} {X_M}{X_M} $$
such that $g|_M \circ \varphi_M = \varphi_M \circ f_M$.
Denote by $\tilde{X}_M \rightarrow X_M$ the normalisation and by \holom{\tilde{\varphi}_M}{\tilde{X}_M}{M}
the induced map, 
then $f_M$ lifts to an endomorphism
\holom{\tilde{f}_M}{\tilde{X}_M}{\tilde{X}_M} such that $g|_{\tilde{M}} \circ \tilde{\varphi}_M = \tilde{\varphi}_M \circ \tilde{f}_M$. 
The general fibre of $\tilde{\varphi}_M$ is a $\PP^1$ or a disjoint union of two $\PP^1$'s and we denote
by the same letter \holom{\tilde{\varphi}_M}{\tilde{X}_M}{\tilde{M}'} its Stein factorisation. It is well-known that this
new fibration is a $\PP^1$-bundle and $\tilde{M}' \rightarrow \tilde{M}$ ramifies exactly in the points where the $\varphi$-fibre 
is a double line. Since such a point gives a singular point of the discriminant locus $\Delta \subset S$ and $\Delta$
is necessarily contained in $\psi$-fibres, Kodaira's classification shows that there are at most two ramification points.
Thus we see that $\tilde{M}'$ is elliptic or rational.  

The $\PP^1$-bundle $\tilde{X}_M \rightarrow \tilde{M}'$ admits an 
endomorphism of degree at least two such that the induced endomorphism on $\tilde{M}'$ is an automorphism
(here we use that $g$ is an automorphism). Thus by Amerik's theorem  \ref{propositionamerik} there exists an \'etale cover
$E \rightarrow \tilde{M'}$ such that the fibre product $\tilde{X}_M \times_{\tilde{M}'} E$ is a product $E \times \PP^1$ and
there exists a finite group $G$ acting diagonally on $E \times \PP^1$ such that $\tilde{X}_M = (E \times \PP^1)/G$.
Thus the projection on the second factor $E \times \PP^1 \rightarrow \PP^1$ induces a fibration
$r_M: \tilde{X}_M \rightarrow \PP^1/G \simeq \PP^1$. 
This fibration descends to a fibration $X_M \rightarrow \PP^1$: this is clear in the complement of the non-normal locus
and the non-normal locus gives a section of $\tilde{X}_M \rightarrow \tilde{M}'$, so by construction it
is contracted by $r_M$.

The construction of the algebraic reduction $r$ is now obvious: let $Y$ be the normalisation of the unique irreducible component of
$\chow{X}$ parametrising a general fibre of the algebraic reduction (the algebraic reduction is almost holomorphic
by Thm. \ref{theoremcampana}, so this is well-defined). Let $\Gamma \rightarrow Y$ be the universal family over $Y$
and denote by $\holom{p}{\Gamma}{X}$ the natural map. Then $p$ is an isomorphism: it is sufficient to check that the
fibre over every point $x \in X$ is a singleton, but this is obvious since $x$ is contained in some $X_M \rightarrow \PP^1$ 
which realises $\Gamma \rightarrow Y$  ``locally''. This proves the {\it claim.}

Since the endomorphism $f$
preserves the algebraic reduction, we clearly have a meromorphic map \merom{g_Y}{Y}{Y}.  
We have constructed the variety $Y$ as the normalisation of some component of $\chow{X}$, thus we can identify
$g_Y$ generically to the push-forward of cycles $f_*$. Since $f_*$ is a holomorphic map on the Chow scheme, we
see that $g_Y$ extends to a holomorphic map. Since $g$ is an automorphism, the endomorphism $g_Y$ has degree at least two.
By \cite{CP00}, there exists a fibration \holom{\tau}{Y}{C}, and
it is easy to check that $g_C \circ \tau = \tau \circ g_Y$. Since $g_C$ is an automorphism and $Y$ is uniruled, 
we deduce by \cite[Lemma 6.1.1,(2)]{Nak08} that $Y$ is a $\PP^1$-bundle over $C$. 
Moreover  another application of Amerik's theorem shows that $Y$ (after \'etale base change)
a product $C \times \PP^1$
Thus the fibre product $Y \times_C S$ is a product $S \times \PP^1$
and the endomorphisms $g$ and $g_Y$ induce an endomorphism $\holom{f'}{S \times \PP^1}{S \times \PP^1}$ of degree at least two.
Since the space of endomorphisms of $\PP^1$ is affine, we have
$f'=(g,h)$ where $\holom{h}{\PP^1}{\PP^1}$  is an endomorphism of degree at least two. 
By the universal property of the fibre product there exists a holomorphic map
$\holom{\mu}{X}{S \times \PP^1}$. We claim that $\mu$ has degree one. 
Since 
$$
b_2(X) = b_2(S)+1 = b_2(S \times \PP^1)
$$
it is obviously finite and therefore $X \simeq S \times \PP^1$.

{\em Proof of the claim.} By construction we have a commutative diagram
$$
 \xymatrix{ 
X \ar[rr]^\mu \ar[d]_f  & & S \times \PP^1 \ar[d]^{f'}
\\
X \ar[rr]^\mu \ar[rd]_\varphi & & S \times \PP^1 \ar[ld]^{p_S}
\\
& S &
}
$$
and we denote by $M$ the branch locus of $\mu$. A straightforward adaptation of Proposition \ref{propositionspecialfibres}
to the case of finite maps shows that $\fibre{f'}{M}=M$. Since $f=(g,h)$ with $h$ of degree at least two,
the divisor $M$ is of the form $\fibre{p_S}{D}$ for some divisor $D$ on $S$.
In particular if $s \in S$ is general point, the restricted map
$$
\holom{\mu|_{\fibre{\varphi}{s}}}{\fibre{\varphi}{s}}{s \times \PP^1}
$$
is not ramified. Since $\PP^1$ is simply connected and $\fibre{\varphi}{s}$ irreducible, the  
restriction $\mu|_{\fibre{\varphi}{s}}$ has degree one. This proves the {\it claim}.

{\it b) $X$ admits a birational contraction $\psi: X \rightarrow X'$.}
Denote by $E$ the exceptional divisor. 
Since the fibers of $X \rightarrow C$ are algebraic, the divisor $E$ must be contained in some fiber $X_c$. 
Hence - at least after passing to some $f^k$ - we have $f^{-1}(E) = E$ and $E$ is not contained in 
the branch locus by Lemma \ref{lemmafibrationcurve}.
Thus by Corollary \ref{corollarycontractionuniruled} the contraction
$\psi$ is the blow-up of a smooth curve $C' \subset X'$ and $X'$ is a compact K\"ahler manifold.
Thus we can proceed inductively and get
a sequence 
$$
X=X_0 \rightarrow X_1 \rightarrow \ldots \rightarrow X_n
$$
such that $X_n$ is compact K\"ahler, admits an endomorphism  \holom{f_n}{X_n}{X_n} of degree $d$ and a fibre type contraction.
By Case a) we know that $X_n \simeq S_n \times \PP^1$.  
We can now argue in the proof of Theorem  \ref{theoremone}
to see that $X \simeq S \times \PP^1$.
\end{proof}

\end{document}